\documentclass[english,11pt,fleqn]{article}
\usepackage{a4wide}
\usepackage{amsmath,amssymb,graphicx}
\usepackage{mathrsfs}
 \usepackage{mathtools}
 \usepackage{float}
\usepackage{dsfont}
\usepackage{comment}
\usepackage[T1]{fontenc}
\usepackage[utf8]{inputenc} 

\def\di{\displaystyle}

\def\div{\mbox{\rm div}}
\newtheorem{theorem}{Theorem}
\newtheorem{definition}{Definition}
\newtheorem{lemma}{Lemma}

\newtheorem{cor}{Corollary}

\newtheorem{remark}{Remark}

\newcommand{\R}{\mathbb{R}}

\newcommand{\C}{\mathbb{C}}

\newcommand{\Cor}{\mbox{\rm Cor}_{\sigma} }
\newcommand{\E}{\mathbb{E}}
\begin{document}
\title{Stochastic modification of Newtonian dynamics and Induced potential - application to spiral galaxies and the dark potential}
\author{Jacky Cresson$^{(a)}$, Laurent Nottale$^{(b)}$, Thierry Lehner$^{(b)}$}
\maketitle
(a) Laboratoire de Math\'ematiques Appliqu\'ees de Pau UMR 5142 CNRS, Universit\'e de Pau et des Pays de l'Adour, E2S UPPA, Avenue de l'Universit\'e, BP 1155, 64013 Pau Cedex, France.

(b) Paris-Sciences-Lettres, Laboratoire Univers et Théories (LUTH) UMR CNRS 8102, Observatoire de Paris, 5 Place Jules Janssen, 92190 Meudon, France.
\begin{abstract}
Using the formalism of stochastic embedding developed by [J. Cresson, D. Darses, J. Math. Phys. 48, 072703 (2007)], we study how the dynamics of the classical Newton equation for a force deriving from a potential is deformed under the assumption that this equation can admit stochastic processes as solutions. We focus on two definitions of a stochastic Newton's equation called differential and variational. 
We first prove a stochastic virial theorem which is a natural generalization of the classical case. The stochasticity modifies the virial relation by adding a potential term called the induced potential which corresponds in quantum mechanics to the Bohm potential. Moreover, the differential stochastic Newton equation naturally provides an action functional which satisfies a stochastic Hamilton-Jacobi equation. The real part of this equation corresponds to the classical Hamilton-Jacobi equation with an extra potential term corresponding to the induced potential already observed in the stochastic virial theorem. The induced potential has an explicit form depending on the density of the stochastic processes solutions of the stochastic Newton equation. It is proved that this density satisfies a nonlinear Schrödinger equation.  Applying this formalism for the Kepler potential, one proves that the induced potential coincides with the ad-hoc "dark potential" used to recover a flat rotation curve of spiral galaxies. We then discuss the application of the previous formalism in the context of spiral galaxies following the proposal and computations given by [D. Da Rocha and L. Nottale, Chaos, Solitons and Fractals, 16(4):565-595, 2003] where the emergence of the "dark potential" is seen as a consequence of the fractality of space in the context of the Scale relativity theory. 
\end{abstract}

\setcounter{tocdepth}{3}

\section{Introduction}

\subsection{General framework}

Let us consider a given deterministic dynamics obtained via a set of physical laws like the Newton fundamental law of dynamics. In classical mechanics, a basic assumption leading to the laws of motion is that a particle describe a differentiable curve in space-time which naturally induces that dynamics is usually described using differential or partial differential equations or more precisely, we restrict our attention to dynamics which can be described using the classical tools of the differential calculus. Note that this assumption of differentiability is an asymptotic one with respect to a given scale of observation for which such a description seems to be valuable (see \cite{cp} for a discussion of this point). However, by definition a rule becomes a law of nature if one can efficiently compare the results with the reality. Doing so, it proves that the initial rule is in some sense robust. The reality is most of the time not described by the ideal mathematical framework used to derive the law. In particular, one expects that stochasticity must be taken into account as the environment of a given experience is never described completely. In particular, even if the law is given in the initial ideal mathematical framework associated to the classical differential calculus, the robustness of the law can be understood as providing some constraints even in a more general setting including stochasticity. As a consequence, one is leaded to give a meaning to a given law expressed using the classical differential calculus over stochastic processes and to study the corresponding dynamics. This is the {\bf extension problem over stochastic processes}. \\

What are the minimal constraints on a given extension ? A first one is that the ideal dynamics, understood as the one obtained in the ideal mathematical framework, must be contained in the extended model, meaning that over stochastic processes which correspond to differentiable functions the extended model reduces to the ideal one. This program is reminiscent of D. Mumford \cite{mumford} call for an integration of stochasticity into the foundations of modeling. It must be noted that the previous problem is very common in many different physical situations as for example in fluid mechanics. The Navier-Stokes equations are derived assuming that the speed of a particle is two times differentiable with respect to space and one time with respect to time. However, real fluids behavior, in particular turbulent ones, are not satisfying these assumptions. As a consequence, one is leaded to extend the notion of solutions by extending the meaning of the equation over a bigger class of objects (Schwartz's distribution in this case). \\

Taking as an extension domain the set of stochastic processes, one can imagine different ways to handle the extension problem.\\

The most common one is to use the {\bf theory of stochastic differential equations} based on the It\^o or Stratonovich stochastic calculus \cite{ok}, considering that the stochasticity induces a {\bf perturbation} term which can be in some cases modeled by a Brownian motion (see \cite{ok}), i.e. informally looking for an equation of the form 
\begin{equation}
\mathscr{P} (d/dt) [x] (t) + "noise"=0, 
\end{equation}
where $\mathscr{P} (d/dt )$ is a differential operator depending on the operator $d/dt$ and "noise" has to be modeled. A typical example in stellar dynamics is given by the seminal work of S. Chandrasekhar (see \cite{chan1}, Chap. IV, Section 1, p.68-70 and \cite{chan2}) where assuming the existence of close encounters between stars in a galaxy the long term motion of a star is modeled by a stochastic differential equation in the Itô sense, i.e. of the form 
\begin{equation}
dX =v(X,t) dt +\sigma (X,t) dW_t ,
\end{equation}
where $W_t$ is a Brownian motion. \\

However, this approach is perturbative in nature and does not give a satisfying answer to the extension problem. Indeed, the new equation is radically different from the initial deterministic one and can not be reduced to it looking for regular solutions.\\

Another strategy is to directly extend the meaning of the dynamics over stochastic processes {\bf without} adding a {\bf perturbation} term but looking for the constraints induced on the stochastic dynamics by the underlying law, i.e. informally by defining a convenient analogue $D$ of $d/dt$ over stochastic processes and looking for 
\begin{equation}
\mathscr{P} (D)[X] =0 ,
\end{equation}
where $X$ is a stochastic process.\\

Such a program was initiated by E. Nelson in \cite{nelson} leading to the framework of {\bf stochastic mechanics}. This framework is an example of {\bf embedding formalism} of differential equations defined in \cite{cr,cr1} in a general setting and in \cite{cd} in the stochastic case. In this article, we follow this strategy. \\

The first step in such a program is to define a convenient framework enable us to write the dynamics. This framework must coincide with the classical differential calculus when the underlying functional space reduces to the set of differentiable functions. Such a formalism is developed in \cite{cd} where a notion of {\bf stochastic derivative} is defined combining the forward and backward derivatives on stochastic processes of E. Nelson \cite{nelson} and using the idea of a complex derivative operator defined by L. Nottale (see \cite{nottale1},p.147 and \cite{nottale2},p.105-106). Such a framework can then be used to extend classical differential equations in the stochastic case using different strategies. Two of them are listed in \cite{cd,cr1}: the differential embedding is a pure algebraic extension obtained by replacing the ordinary time derivative by the stochastic derivative. The variational embedding consists in extending first the notion of Lagrangian functional and then developing the corresponding stochastic calculus of variations. Having such an extension of a given dynamics over stochastic processes, one can then study the persistence of some properties of the initial deterministic dynamics like symmetries, first integrals, etc. The validity of such an extension in order to describe the effective dynamics of a given physical process must be cases by cases discussed.\\ 

The stochasticity can be intrinsic like in some views of quantum mechanics \cite{nelson} or extrinsic like a convenient description of the long term behavior of chaotic dynamical systems \cite{nottale3,cr}. In some applications, the stochastic embedding formalism (differential or variational) can be seen as a concrete realization of the {\bf scale relativity} principle of L. Nottale \cite{nottale1}: equations of physics must keep the same form at each scale, the word "same" being understood as keeping the algebraic form of the differential operator (differential embedding) or the variational character (variational embedding) or both (coherent embedding) (see \cite{cr1,cd}).

\subsection{Main results}

In this article, we focus on the stochastic embedding of the Newton equation for a given potential $U$. We first prove that the solutions of the (differential) stochastic Newton equation are always gradient diffusion thanks to the characterization of gradient diffusion given by S. Darses and I. Nourdin in \cite{dn1}. This result induces a natural complex valued action functional which generalize the classical action functional in mechanics. In particular, we write a stochastic version of the Hamilton-Jacobi equation. The real part of this equation can be viewed as a modified Hamilton-Jacobi equation the modification being on the potential which is in some sense corrected by a stochastic "induced" potential. The imaginary part corresponds to a continuity equation satisfies by the density of the stochastic process solution of the stochastic Newton equation. Many applications in mechanics make use of the virial Theorem relating the kinetic and potential energy at equilibrium. A typical application in Astrophysics is the determination of the mass of a galaxy from the observation of the characteristic speed of stars. We derive stochastic analogues of the virial Theorem depending on the notion of equilibrium one is looking for. In the strong case, we look for a quantity which is a constant stochastic process. This leads to a result similar to the classical virial Theorem but again the potential is corrected by the same potential as for the stochastic Hamilton-Jacobi equation of pure stochastic origin. The second point of view is to look for quantities whose expectation are constant called the weak case. Another result is then obtained but more difficult to interpret. The stochastic induced potential depends explicitly on the density of the stochastic process solution of the stochastic Newton equation. In order to compute this density, we make use of the result obtained by one of the author and S. Darses in \cite{cd} relating the stochastic Newton equation to the Schr\"odinger equation. This connection allows us to use known results on the solutions of the Schr\"odinger equation to obtain information of the density of the stochastic processes solution of the stochastic Newton equation. This take this opportunity to cancel the initial assumption of gradient diffusion made in \cite{cd} and going back originally to an assumption of E. Nelson \cite{nelson2} in the context of stochastic mechanics, as it follows from the structure of the equation. We then apply our results taking the Kepler potential. In that case, we give an explicit form to the stochastic induced potential and we derive several dynamical results. In particular, the flat rotation curve Theorem which shows that the real part of the speed is constant in that case. This constant speed can be used to obtain an explicit form of the diffusion coefficient. Our main objective being to apply our result in the context of spiral galaxies, we derive a Noether type result which enables us to obtain an estimation of the circular velocity of a particle moving on a circular orbit. We then apply all these results in the context of the dynamics of spiral galaxies. Using the Scale Relativity Theory developed by L. Nottale in \cite{nottale1,nottale2} as a support to justify the use of the stochastic Newton equation and using known numerical values for mass and speed of the observed rotation curve for the Milky way, we predict the distance at which the rotation curve begins to be flat and obtain a very good agreement with the observations. Moreover, the stochastic induced potential corresponds in that case to the ad-hoc "dark potential" used in the literature to recover the observed rotation curve of spiral galaxies.

\subsection{Connection with previous results}

The stochastic embedding formalism used in this article use the tools developed by E. Nelson for the formulation of the stochastic mechanics \cite{nelson, nelson2,nelson3}. The main difference is the use of the stochastic derivative mixing in a complex valued operator the left and right stochastic derivative of Nelson following an idea of L. Nottale in \cite{nottale1}. We refer to \cite{cd} where this formalism is developed. Several computations and results are derived in an informal way in different articles and books of L. Nottale \cite{nottale1,nottale2}. In particular, the stochastic embedding formalism can be seen as a formal account of the strategy described by L. Nottale in \cite{nottale1} to discuss the connection between the scale relativity theory and stochastic mechanics and in \cite{nottale3} relating the long term dynamics of dynamical systems with a stochastic system.\\

The strong stochastic virial Theorem was discussed in the context of stochastic mechanics by S.M. Moore in \cite{moore} but the proof seems to be incomplete (see Section \ref{virial}. The weak stochastic virial Theorem is discused by P.H. Chavanis in (\cite{chavanis},V) in a particular case. \\

A real stochastic Hamilton-Jacobi equation was discussed by F. Guerra and L.M. Morato \cite{guerramorato1} in the context of stochastic mechanics. They obtain this equation from a particular notion of critical point for a functional called in the sense of Lafferty by E. Nelson (see \cite{nelson3},p.439). The main difference is the complex nature of our stochastic Hamilton Jacobi equation which implies that we obtain some mixing between the real and imaginary part of the complex derivative in the equation. Another one (also real) was introduced by J-P. Ortega and L. Cami in \cite{cami} in the context of stochastic geometric mechanics in \cite{cami} following the seminal work of J-M. Bismut \cite{bismut} on random mechanics. The two results do not coincide as the underlying theories are completely different.\\ 

The stochastic induced potential appears already in the literature and is known as the Bohm's potential or quantum potential. We refer to the book of J-C. Zambrini (\cite{zambrini2},p.168) for an historical account of the emergence of this potential in quantum mechanics related with the hidden variable theory of D. Bohm. \\

The connection between the stochastic Newton equation used in this article and the Schr\"odinger equation was already discussed in \cite{cd} and is based on the seminal work of E. Nelson \cite{nelson} and the idea of L. Nottale in \cite{nottale1}. Other derivation of the same result are done for example by P-H. Chavanis in (\cite{chavanis}, equation (B.1)). These derivations make the assumption that the stochastic derivative is a gradient or give an informal proof. In this article, we cancel the gradient assumption as it follows from the structure of the equation. \\
 
The application to the dynamics of spiral galaxies follows the previous work of L. Nottale and D. Rocha in \cite{rocha} and (\cite{nottale2},Section 13.8.2 p.652-654). All their results are analyzed in the stochastic framework that we develop. \\

A very close related work is done by F. Pierret in \cite{pierret} using the scale dynamics formalism developed in \cite{cp}. The framework is not stochastic. A generalized version of the virial Theorem as well as a generalized Hamilton-Jacobi equation are derived. In particular, the stochastic induced potential appears in \cite{pierret} as a scale dynamical effect. 

\subsection{Organization of the paper}

The plan of the paper is as follows: \\

In Section \ref{remind}, we remind the definition of the stochastic derivative introduced in \cite{cd} and some of its properties. In Section \ref{sectionstocnewton}, we define the differential and variational stochastic Newton equation following \cite{cd} and the stochastic embedding formalism. In Section \ref{virial}, we prove a stochastic version of the virial Theorem. Section \ref{properties} is dedicated to some properties of the differential stochastic Newton equation. In particular, we prove that in this case, solutions of the equation are gradient diffusions. Section \ref{sectionaction}, we define a stochastic action functional and we derive the corresponding form of the stochastic Hamilton-Jacobi equation. Section \ref{sectionemergence} shows that the real part of the stochastic Hamilton-Jacobi equation is a perturbation of the classical Hamilton-Jacobi equation by a potential which is induced by the stochastic character of the underlying set of solutions. This induced potential can be explicitly written in term of the density of the stochastic process solution of the stochastic Newton equation. In order to compute the density and to explicit the induced potential, we transform in Section \ref{sectionschro} the stochastic Newton equation in a Schr\"odinger like equation using a change of variable. As an application, we explicit the form of the induced potential when the underlying potential is the Kepler one. We then obtain that the induced potential corresponds to the ad-hoc choice of a "dark potential" used in the literature to explain the flat rotation curve of spiral galaxies. The application of such a formalism in this setting is discussed. Finally, we give in Section \ref{sectionpers} some perspectives of this work.

\section{Reminder about the stochastic derivative}
\label{remind}

In this Section, we remind some definitions and properties obtained in \cite{nelson,cd,darses2,dn1} about Nelson left and right derivatives and the {\bf stochastic derivative} introduced in \cite{cd}. We do not give the more general class of stochastic processes for which such quantities are defined and for which formula can be given but instead restrict our attention to Brownian diffusion processes with a constant diffusion. We refer to \cite{cd,darses2,dn1} for more general results.

\subsection{The stochastic derivative}

In the following, we restrict our attention to stochastic processes which are diffusion processes of the form 
\begin{equation}
\label{browndiff}
dX_t =v(t,X_t) dt +\sigma dW_t ,
\end{equation}
where $W_t$ is a $\R^n$ Brownian motion, $\sigma \in \R$ is constant and $v:\R \times \R^n \rightarrow \R^n$ is a function such that 
\begin{itemize}
\item (i) For all $x,y\in \R^d$, $\sup_{t\in [0,T]} \mid v(t,x)-v(t,y)\mid \leq K \mid x-y\mid$, 

\item (ii) $\sup_{t\in [0,T]} \mid v(t,x) \mid \leq K (1+\mid x \mid )$,
\end{itemize}
\vskip 2mm
which are classical conditions (see \cite{ok} Theorem 5.2.1 p.70) ensuring the existence and uniqueness of solutions for It\^o stochastic differential equations.

\begin{remark}
The terminology of {\bf diffusion coefficient} is often used to speak of $\sigma$. This is due to the fact that if one considers the Brownian motion $dX_t =\sigma dW_t$, the density $p_t (x) =u(t,x)$ of $X_t$ satisfies the {\bf diffusion equation} or {\bf heat equation} 
\begin{equation}
\di\frac{\partial u}{\partial t } = \di\mathscr{D} \Delta u  ,
\end{equation}
with diffusion coefficient
\begin{equation}
\label{relnotdif}
\mathscr{D} = \di\frac{\sigma^2}{2} .
\end{equation}
We refer to the historical account given in (\cite{evans},Section 3.1.1 p.37-39) for more details. \\

The coefficient $\mathscr{D}$ is used by L. Nottale in (\cite{nottale1}, p.143-153) to discuss the connection between stochastic mechanics and scale relativity. However, the diffusion coefficient $\mathscr{D}$ has a completely different meaning. It is understood in scale relativity as a parameter related to the {\bf fractal dimension} (see \cite{tricot}) of the trajectories. As a consequence, it is related to a fundamental geometric property of space-time which is assumed to be fractal in scale relativity.\\

It must be noted that the assumption that trajectories have a fractal dimension equal to $2$ is more general than the assumption that the trajectories are described by stochastic diffusion processes. In particular, all the geometric properties of trajectories in scale relativity can not be reducible to the use of stochastic processes, even by extending the class of these processes in order to include fractional Brownian motion of arbitrary order. We refer to (\cite{nottale1},p.145-146) for more details.
\end{remark}

We denote by $D_+$ and $D_-$ the right and left Nelson derivatives (see \cite{nelson}):

\begin{equation}
D_+ [X_t ]=\di\lim_{h\rightarrow 0^+} \mathbb{E} \left [ 
\di\frac{X_{t+h} -X_t}{h} \mid \mathcal{P}_t \right ] , \ \ \
D_- [X_t ]=\di\lim_{h\rightarrow 0^+} \mathbb{E} \left [ 
\di\frac{X_{t} -X_{t-h}}{h} \mid \mathcal{F}_t \right ] ,
\end{equation}
where $\mathcal{P}_t$ (resp. $\mathcal{F}_t$) is the forward (resp. backward) differentiating $\sigma$-field for Brownian diffusions $X$ of the form (\ref{browndiff}) in the sense of \cite{darses2}.\\

Under weak conditions, one can prove that  (See \cite{dn1},Proposition 2):
\begin{equation}
D_+ [X_t ] =v(t, X_t ) ,\ \ D_- [X_t ] = v (t,X_t ) -\sigma^2 \di\frac{\nabla p_t}{p_t} (X_t ) ,
\end{equation}
where $p_t$ is the density associated to the law of $X_t$. \\

Moreover, one can prove the following chain rule formula (see \cite{dn1},Proposition 3):
\begin{equation}
\label{comp}
D_{\pm} [g(X_t ,t) ] =\di\left (
\partial_t +D_{\pm} [X_t] \, \nabla_x \pm \di\frac{\sigma^2}{2} \Delta_x \right ) 
[g] (X_t ,t) ,
\end{equation}

The {\bf stochastic derivative} introduced in \cite{cd} is defined as
\begin{equation}
\mathscr{D}_{\mu} =\di\frac{D_+ +D_-}{2} +i\mu \frac{ D_+ - D_-}{2} , \mu \in \{ -1 ,1 \} .
\end{equation} 
This operator is extended to complex stochastic process by linearity. This operator appears first in the work of L. Nottale (see [\cite{nottale1},p.147] and [\cite{nottale2},p.165-166]).

\begin{remark}
The complex nature of the stochastic derivative, i.e. the special form of the new operator combining the two information contained in $D_+ X$ and $D_- X$, is supported both by mathematical arguments related to the notion of doubling algebra (see \cite{nottale2},p.160-161) and physical arguments (see \cite{nottale2},p.161-164) related to the requirement of form invariance of fundamental equations or covariance as used by A. Einstein \cite{ein}.
\end{remark}

Two comments on the form of the stochastic derivative:\\

\begin{itemize}
\item The stochastic derivative gives the same weight to $D_+$ and $D_-$, i.e. to the past and the future.

\item If $X$ is a deterministic differentiable process, i.e. $dX_t (\omega )=f(t) dt$ for all $\omega \in \Omega$, then $D_+ X =D_- X$ and $\mathscr{D}_{\mu} X = f(t)$, which coincides with $\di\frac{d X_t}{dt}$. 
\end{itemize}

As a consequence, the stochastic derivative can be used for the extension problem stated in the Introduction.\\

The set of diffusion processes for which the stochastic derivative is well defined coincides with the one ensuring the existence of the Nelson backward and forward derivatives $D_+$ and $D_-$. This set has been well studied in the literature and we refer to \cite{cd} for an overview and in particular to the work of E. Nelson \cite{nelson}, W.A. Zheng and P.A. Meyer \cite{mayer1,mayer2}, R. Carmona \cite{carmona} and H. F\"ollmer \cite{follmer} and more recently to the work of L. Wu \cite{wu} and S. Darses and I. Nourdin \cite{dn1,darses2}.

\subsection{Properties of the stochastic derivative}

Using the previous results about Nelson's derivatives, we can prove that 
\begin{equation}
\label{correction}
\mathscr{D}_{-\mu} (X_t ) =\mathscr{D}_{\mu} (X_t ) +i\mu \Cor (X_t ) ,
\end{equation}
where the correction term $\Cor (X_t )$ is given by 
\begin{equation}
\Cor (X_t ) = \sigma^2 \di\frac{\nabla p}{p} .
\end{equation}
We also obtain the following chain rule formula
\begin{equation}
\label{chaine}
\mathscr{D}_{\mu} [g(X_t ,t )] = \left ( \partial_t +\mathscr{D}_{\mu} (X_t ) \nabla +i\mu \di\frac{\sigma^2}{2} \Delta \right )  [g] (X_t ,t ) ,
\end{equation}
which leads to 
\begin{equation}
\mathscr{D}_{-\mu} [g(X_t ,t)]= \mathscr{D}_{\mu} [g(X_t ,t)] +i\mu \left ( \Cor (X_t ) \nabla -\sigma^2 \Delta \right ) [g](X_t ,t).
\end{equation}
The two expressions can be encoded in a single formula as 
\begin{equation}
\label{chainecor}
\mathscr{D}_{\alpha\mu} [g(X_t ,t)]= \mathscr{D}_{\mu} [g(X_t ,t)] +i\di\frac{(1-\alpha )}{2} \mu \left ( \Cor (X_t ) \nabla -\sigma^2 \Delta \right ) [g](X_t ,t),
\end{equation}
where $\alpha =\pm 1$. \\

A Leibniz like formula can be obtained for the stochastic derivative: let $X$ and $Y$ be two complex stochastic processes, then 
\begin{equation}
\label{leibniz}
\E \left [ \mathscr{D}_{\mu} X \cdot Y +X\cdot \mathscr{D}_{-\mu} Y \right ] =\di\frac{d}{dt} \E (X\cdot Y ) ,
\end{equation}
where $\E$ is the expectation. The interplay between $\mathscr{D}_{\mu}$ and $\mathscr{D}_{-\mu}$ comes from the non reversibility between $D$ and $D_*$.\\

In particular, replacing $\mathscr{D}_{-\mu} Y$ by its expression (\ref{correction}) in function of $\mathscr{D}_{\mu}$, we obtain 
 \begin{equation}
\label{leibnizcor}
\E \left [ \mathscr{D}_{\mu} X \cdot Y +X\cdot \mathscr{D}_{\mu} Y \right ] +
i\mu \E \left [ X \cdot \Cor (Y) \right ] =\di\frac{d}{dt} \E (X\cdot Y ) ,
\end{equation}

Another important property that we will use in the following is the composition lemma:

\begin{lemma}[Composition] 
For $\alpha =\pm 1$ and $\mu=\pm 1$, we have 
\begin{equation}
\left .
\begin{array}{lll}
\mathscr{D}_{\alpha \mu} \circ \mathscr{D}_{\mu} & = &  
\di\frac{1}{4} \left [ (D_+^2 +D_-^2 ) (1-\alpha ) +(D_+ D_- +D_- D_+ ) (1+\alpha ) \right ] \\
 & & +i \mu \di\frac{1}{4} 
\left [ (D_+^2 -D_-^2 ) (1+\alpha ) +(D_- D_+ -D_+ D_- ) (1-\alpha ) \right ] .
\end{array}
\right .
\end{equation}
\end{lemma}

When $\alpha=1$ we obtain for the real part of $\mathscr{D}_{\mu}^2$ the quantity $\di\frac{1}{2} (D_+ D_- +D_- D_+ )$ called {\bf mean second derivative} or {\bf mean acceleration} by E. Nelson in (\cite{nelson}, equation (15) p.99). The imaginary part reduces to $\di\frac{\mu}{2} (D_+^2 -D_-^2 )$ which will play an essential role in the construction of an action functional in Section \ref{sectionaction}.

\section{The differential and variational stochastic Newton equation}
\label{sectionstocnewton}

In this Section, we discuss different versions of what can be called a stochastic Newton equation. All these equations can be obtained via the stochastic embedding of dynamical systems defined in \cite{cd} and differ in the properties they are preserving from the classical Newton equation. 

\subsection{The Newton equation}

The Newton equation also called fundamental equation of dynamics is the differential equation
\begin{equation}
\di\frac{d}{dt} \left ( \partial_v K (\dot{x} )\right )=F(x ,\dot{x}), 
\end{equation}
where $\dot{x}=\di\frac{dx}{dt}$, $K$ is a homogeneous function of order $2$ and $F$ is a force which can depends on $x$ and linearly on $\dot{x}$. \\

Typical examples of forces are given by a force deriving from a potential $U :\R^n \rightarrow \R$ which is a given $C^1$ function such that 
\begin{equation}
F(x,\dot{x} ) =-\nabla U ,
\end{equation}
and friction forces which are of the form 
\begin{equation}
F(x,\dot{x})=-\gamma \dot{x} .
\end{equation}

When the force is conservative, meaning that the force derives from a potential $U$, one can prove that the solutions of the Newton equation are in correspondence with critical points of a functional defined by 
\begin{equation}
\label{lagrangefunc}
\mathscr{L} (x)=\di\int_a^b L (x (s),\dot{x} (s) ) \, ds,
\end{equation}
where $L :\R^n \times \R^n \rightarrow \R$ is called a {\bf natural Lagrangian} and is given by
\begin{equation}
L(x,v)=K(v)-U(x) .
\end{equation}
A stochastic extension of the Newton equation must then include a discussion of the properties the new equation is preserving. We discuss this point in the following Section using the formalism of the stochastic embedding developed in \cite{cd}.

\subsection{The differential stochastic Newton equation}
\label{newtondiff}

A first idea to generalize the Newton equation is to preserve its form, meaning the algebraic structure of the differential operator by replacing the classical derivative by our stochastic derivative. Such a procedure is called a differential embedding in \cite{cd}. 

We restrict our attention in all this Section to Brownian diffusion
\begin{equation}
\label{proc}
dX_t = v(X_t ,t ) dt +\sigma dW_t ,
\end{equation}
where $\sigma$ is a constant and $W_t$ is a Brownian motion. \\

The differential embedding of the Newton equation then leads to 
\begin{equation}
\label{diffnewton}
\mathscr{D}_{\mu} \left ( \partial_v K (\mathscr{D}_{\mu} X_t ) \right )=F (X_t ,\mathscr{D}_{\mu} X_t ) ,\ \mu\in \{ -1,1\} . 
\end{equation}
A main remark is that due to the complex nature of the stochastic derivative the left hand-side is a priori complex. The nature of the force has then huge consequences. In particular, we have the following distinctions:\\

\begin{itemize}
\item A force deriving from a potential $U (x)$ will always produce a real quantity by a stochastic embedding as the stochastic embedding of $F(x)=-\nabla U (x)$ is given by $F(X_t )=-\nabla U$ which is again real. 

\item A dissipative force, depending on the speed of $x$, as for example a friction or damping term $F(x)=-\gamma \dot{x}$ will produce by the stochastic embedding a term given by $F(X_t ) =-\gamma \mathscr{D}_{\mu} X_t$ which is in general a complex quantity.
\end{itemize}

As we will see in Section \ref{properties}, this property will cancel the possibility to define a natural action functional in the dissipative case.

\subsection{The variational stochastic Newton equation}

We restrict our attention to the conservative version of the Newton equation. The stochastic embedding of the functional $\mathscr{L}$ defined in (\ref{lagrangefunc}) is defined by 
\begin{equation}
\label{stocfunc}
\mathscr{L} (X_t ) =\di\E \di \left ( \int_a^b L(X_s ,\mathscr{D}_{\mu} X_t ,s )\, ds \right ) ,
\end{equation}
which is a generalization of a seminal work of K. Yasue \cite{yasue} where a variational formulation is obtained for the stochastic Newton equation defined by E. Nelson \cite{nelson}. \\

Developing the stochastic calculus of variations due to K. Yasue \cite{yasue} in our setting, we have proved that the stochastic Euler-Lagrange associated to the functional (\ref{stocfunc}) when $L$ is a natural Lagrangian is given by 
\begin{equation}
\label{variationalnewton}
\mathscr{D}_{-\mu} \left ( \partial_v K (\mathscr{D}_{\mu} X_t ) \right )=-\nabla (X_t) ,\ \mu\in \{ -1,1\} ,
\end{equation}
and called the variational stochastic Newton equation.\\

As we see, the variational stochastic Newton equation does not correspond to the differential stochastic Newton equation. The dependence of the equation on the extension one is performing is called the coherence problem and is in fact a general problem in all embedding formalism \cite{cr1}.

\subsection{Combining the two forms of stochastic Newton equations}

The previous definitions can be mixed in a single one by introducing a constant $\alpha =\pm$. The $\alpha$-stochastic Newton equation 
\begin{equation}
\label{dynsto}
\mathscr{D}_{\alpha \mu} \left ( \partial_v K (\mathscr{D}_{\mu} X_t ) \right )=F (X_t ,\mathscr{D} X_t ) ,\ \mu\in \{ -1,1\}, 
\end{equation}
and $\alpha =\pm 1$.\\

When $\alpha=1$ this equation corresponds to the stochastic differential embedding of the classical Newton equation. When $\alpha=-1$ this equation corresponds to the stochastic variational embedding of the Newton equation and we have a one to one correspondence between critical points of a Lagrangian functional and solutions of the $-1$-stochastic Newton equation.

\section{A stochastic virial theorem}
\label{virial}

In this Section, we derive a stochastic version of the classical virial Theorem. However, due to the existence of different notions of first integrals in the stochastic case (strong or weak), one can look for different type of generalization. This point is discussed in the following section \label{weakstrong}. We then derive the strong stochastic virial Theorem whose interpretation is clear and the weak one whose dynamical meaning is not simple. Previous results in these directions have been obtained by S.M. Moore in \cite{moore} where a strong virial theorem is discussed in the context of stochastic field theory using the stochastic mechanics of E. Nelson \cite{nelson} and by P.H. Chavanis in (\cite{chavanis},V) where weak virial theorem is discused in a particular case.

\subsection{Strong versus weak stochastic virial Theorem}

Let $K$ be the homogeneous function of order $2$ given by 
\begin{equation}
K(v)=\di\frac{1}{2} m v^2 ,
\end{equation}
and $U$ be a homogeneous function of order $\gamma \in \R$. The classical virial theorem asserts that 
\begin{equation}
\label{virialequa}
\di\frac{d}{dt} (m x. \dot{x} ) =2K - \gamma U ,
\end{equation}
so that when the system is at equilibrium we obtain 
\begin{equation}
2K =\gamma U .
\end{equation}
Another formulation of the same result can be obtained remarking that the moment of inertia $I(x)=m x^2$ satisfies 
\begin{equation}
\label{virialequa2}
\di\frac{1}{2} \di\frac{d^2 I}{dt^2} = 2K-\gamma U ,
\end{equation}
and that at equilibrium one has 
\begin{equation}
\di\frac{d^2 I}{dt^2} =0 , 
\end{equation}
meaning that $\di\frac{dI}{dt}$ is a first integral of the system.\\ 

Looking for a stochastic version of the virial Theorem several possibilities can be considered depending on the object one focus and the way we interpret equilibrium.\\

Let $J(x , v)$ be a quantity depending on $x$ and $v$ which in the classical case is evaluated over solutions of the system as $J(x,\dot{x})$. What is the stochastic analogue of the quantity 
\begin{equation}
\label{virial3}
\di\frac{d}{dt} [ J(x,\dot{x} )]\ \ ? 
\end{equation}

First, the stochastic embedding of $J$ leads directly to the following stochastic quantity $J(X_t , \mathscr{D}_{\mu} X_t )$ over the solution of the stochastic Newton equation. Second, generalizing equation (\ref{virial3}) one can look for 
\begin{equation}
\label{virial4}
\mathscr{D}_{\mu} \left [ J( X_t , \mathscr{D}_{\mu} X_t ) \right ] , 
\end{equation}
or to an averaged quantity
\begin{equation}
\label{virial5}
\di\frac{d}{dt} \left [ \mathbb{E} \left ( J( X_t , \mathscr{D}_{\mu} X_t ) \right ) \right ] .
\end{equation}

In the next Sections, we explore these two possibilities.

\subsection{Strong stochastic virial Theorem}

The following result, which is a consequence of (\cite{nelson},Theorem 11.11 p.96), can be useful to choose a possibility:

\begin{lemma}
Let $Y$ be a process such that $\mathscr{D}_{\mu} Y$ exists. Then $Y$ is a constant (i.e. $Y_t$ is the same random variable for all $t$) if and only if $\mathscr{D}_{\mu} Y=0$.
\end{lemma}

If we want to keep the classical interpretation of the virial Theorem view as the fact that the derivative of the moment of inertia is a constant of motion, we can compute first $\mathcal{D}^2 (mX^2 )$. 

\begin{theorem}[Strong stochastic virial theorem]
Let $X$ be satisfying the $\alpha$-Newton equation (\ref{dynsto}), then we have 
\begin{equation}
\mathscr{D}_{\mu}^2 (m X^2 )=4K -2\gamma U +2i\mu m\sigma^2 \mbox{div} (\mathscr{D}_{\mu} X ) .
\end{equation}
\end{theorem} 

\begin{remark}
A stochastic virial Theorem is derived by S. M. Moore in (\cite{moore}, Appendix p.2108-2109) which has some connections with our result. However, we have not been able to check the computations of the paper. In particular, equality (66) of \cite{moore} seems a priori non trivial and does not follow from the classical properties of the Nelson forward and backward derivatives.
\end{remark}

\vskip 2mm {\bf Proof}. We denote the two real components of $\mathscr{D}_{\mu} X$ by $\mathbf{v}$ and $\mathbf{u}$, i.e. 
\begin{equation}
\mathscr{D}_{\mu} X =\mathbf{v}+i\mu \mathbf{u} ,
\end{equation}
where $\mathbf{v}$ and $\mathbf{u}$ are functions of $X$ and $t$ and correspond to the {\bf current velocity} and {\bf osmotic velocity} defined by E. Nelson in (\cite{nelson}, p.105-106). \\

Using the chain rule formula (\ref{chaine}), with $g(x)=x^2$, we obtain
\begin{equation}
\mathscr{D}_{\mu} (mX^2 ) =2m (X.\mathscr{D}_{\mu} X ) +im\mu \sigma^2 .
\end{equation}
As a consequence, we have 
\begin{equation}
\mathscr{D}_{\mu}^2 (mX^2 )=2m \mathscr{D}_{\mu} \left ( X. \mathscr{D}_{\mu} X \right ) .
\end{equation}
Using the decomposition of $\mathscr{D}_{\mu} X$ as a function of $\mathbf{v}$ and $\mathbf{u}$, we are lead to the computation of two quantities 
\begin{equation}
\mathscr{D}_{\mu} (X. \mathbf{v} )\ \ \ \mbox{and}\ \ \ \mathscr{D}_{\mu} (X.\mathbf{u} ) .
\end{equation}
Here again we use the chain rule formula with $g_1 (x,t)=x \mathbf{v} (x,t)$ and $g_2 (x,t)=x. \mathbf{u} (x,t)$ respectively. We then obtain 
\begin{equation}
\mathscr{D}_{\mu} (X. \mathbf{v} ) =X\di\frac{\partial \mathbf{v}}{\partial t}+\mathscr{D}_{\mu} X \left ( 
\mathbf{v} +X \di\frac{\partial v}{\partial x} \right ) +i\mu\di\frac{\sigma^2}{2} 
\left ( 2\di\frac{\partial \mathbf{v}}{\partial x} +X\di\frac{\partial^2 \mathbf{v}}{\partial x^2} \right ) , 
\end{equation}
and a similar equation for $\mathscr{D}_{\mu} (X.\mathbf{u} )$. Combining the two equations, we finally obtain 
\begin{equation}
\label{intervirial}
\mathscr{D}_{\mu}^2 (mX^2 )=2m (\mathscr{D}_{\mu} X )^2 +i 2m \mu \mbox{div} \left ( \mathscr{D}_{\mu} X \right ) + 2mX \left ( \mathbb{L} [\mathbf{v} ] +i\mu \mathbb{L} [\mathbf{u} ] \right ) , 
\end{equation}
where $\mbox{div} (\mathscr{D}_{\mu} X ) = \partial_x \mathbf{v} +i\mu \partial_x \mathbf{u}$ and $\mathbb{L}$ is the differential operator
\begin{equation}
\mathbb{L} =\partial_t  +\mathscr{D}_{\mu} X \partial_x +i\mu \di\frac{\sigma^2}{2}\Delta .
\end{equation}
As $\mathscr{D}_{\mu}^2 X =\mathscr{D}_{\mu} ( \mathbf{v} +i\mu \mathbf{u} )$, one shows that 
\begin{equation}
\mathbb{L} [\mathbf{v} ] +i\mu \mathbb{L} [\mathbf{u} ] =\mathscr{D}_{\mu}^2 X .
\end{equation}
Equation (\ref{intervirial}) can then be written as 
\begin{equation}
\mathscr{D}_{\mu}^2 (mX^2 )=2m (\mathscr{D}_{\mu} X )^2 +i 2m \mu \mbox{div} \left ( \mathscr{D}_{\mu} X \right ) + 2mX \mathscr{D}_{\mu}^2 X . 
\end{equation}
As $X$ is a solution of the stochastic Newton equation, we have $m X \mathscr{D}_{\mu}^2 X =-X \nabla U$. Due to the homogeneity of $U$, we then have due to the Euler Theorem that 
\begin{equation}
m X \mathscr{D}_{\mu}^2 X =-\gamma  U .
\end{equation}
As a consequence, using that $K (v)=\di\frac{1}{2} m v^2$, we obtain 
\begin{equation}
\mathscr{D}_{\mu}^2 (mX^2 )=4 K (\mathscr{D}_{\mu} X ) +i 2m \mu \mbox{div} \left ( \mathscr{D}_{\mu} X \right ) -2\gamma U . 
\end{equation}
This concludes the proof. $\square$\vskip 1mm

Assuming that at {\bf equilibrium} one has 
\begin{equation}
\mathscr{D}_{\mu}^2 X =0 ,
\end{equation}
we conclude that 
\begin{equation}
\label{virialstrong}
2 K (\mathscr{D}_{\mu} X ) +i m \mu \mbox{div} \left ( \mathscr{D}_{\mu} X \right ) -\gamma U =0 . 
\end{equation}
The real part of equation (\ref{virialstrong}) is equivalent to 
\begin{equation}
m\mathbf{v}^2 = \gamma  U +m \left [ \mathbf{u}^2 +\sigma^2 \mbox{div} \mathbf{u} \right ] .
\end{equation}
The stochastic component of the dynamics induced an extra potential given by 
\begin{equation}
U_{\sigma ,induced} = -m \left [ \mathbf{u}^2 +\sigma^2 \mbox{div} \mathbf{u} \right ] . 
\end{equation}
Using the fact that $\mathbf{u}$ can be written as a function of the density $p_t$ of the stochastic process $X$ as 
\begin{equation}
\mathbf{u} =\di\frac{\sigma^2}{2} \di\frac{\nabla p_t }{p_t} ,
\end{equation}
we deduce that the induced potential in the strong stochastic virial Theorem is given by 
\begin{equation}
U_{\sigma ,induced} = -m\di\frac{\sigma^4}{2} \di\frac{\Delta (\sqrt{p} )}{\sqrt{p}} .
\end{equation}
We then deduce that the strong stochastic virial Theorem gives at equilibrium the relation 
\begin{equation}
m\mathbf{v}^2 = \gamma  U -U_{\sigma ,induced} .
\end{equation}
The extra potential $U_{\sigma ,induced}$ will reappear in Section \ref{sectionemergence} as an extra potential term in the real part of the stochastic Hamilton-Jacobi equation. 

\subsection{Weak stochastic virial Theorem}

In this Section, we look for a stochastic analogue of the relation (\ref{virialequa}) by computing the quantity
\begin{equation}
\di\frac{d}{dt} 
\left ( 
\E (X\cdot \nabla K ) 
\right ) .
\end{equation} 
We obtain the following result:

\begin{theorem}[Weak Stochastic virial theorem] 
Let $X$ satisfying the $\alpha$-Newton equation (\ref{dynsto}), then we have 
\begin{equation}
\di\frac{d}{dt} 
\left ( 
\E (X\cdot \nabla K ) 
\right ) 
= 
\E 
\left (  
2K -\gamma U +i\mu \di\frac{(1+\alpha)}{2} (\Cor (X_t ) \cdot \nabla K )
\right ) .
\end{equation}
\end{theorem}

As we see, there exists a strong difference between the variational Newton equation corresponding to $\alpha =-1$ where no modifications of the classical virial Theorem exist and the differential Newton equation corresponding to $\alpha=1$ where a correction term emerges. In the variational setting, the fact that the formulation of the virial theorem is preserved must be put in correspondence with the persistence of a Noether type result in the stochastic case. We refer to \cite{cd} for more details.

\vskip 2mm {\bf Proof}. 
We have by the chain rule for Nelson's derivatives that 
\begin{equation}
\E (X.\mathscr{D}_{\alpha \mu} ( \nabla K  ))= \di\frac{d}{dt} \left ( \E (X\cdot \nabla K ) \right ) - \E ( \mathscr{D}_{-\alpha \mu} X \cdot \nabla K ) .
\end{equation}
The homogeneity of $K$ implies that 
\begin{equation}
\mathscr{D}_{\mu} X_t \cdot \nabla K (\mathscr{D}_{\mu} X_t ) = 2 K (\mathscr{D}_{\mu} X_t ) ,
\end{equation}
and the homogeneity of $U$ that 
\begin{equation}
X\cdot \nabla U = \gamma U .
\end{equation}
Multiplying the $\alpha$-stochastic Newton equation by $X$ and taking the expectation, we obtain 
\begin{equation}
\E (X.\mathscr{D}_{\alpha \mu} \nabla K  )
= - \E ( X. \nabla U ) .
\end{equation}
For the variational stochastic Newton equation corresponding to $\alpha =-1$, we obtain using the chain rule formula 
\begin{equation}
\left .
\begin{array}{lll}
\di\frac{d}{dt} \left ( \E (X\cdot \nabla K ) \right ) & = &  \E (\mathscr{D}_{\mu} X \nabla K ) -\E (X\cdot \nabla U ) , \\
 & = & \E (2K -\gamma U ) .
\end{array}
\right .
\end{equation}
For the differential stochastic Newton equation corresponding to $\alpha =1$, we have 
\begin{equation}
\left .
\begin{array}{lll}
\di\frac{d}{dt} \left ( \E (X\cdot \nabla K ) \right ) & = &  \E (\mathscr{D}_{-\mu} X \nabla K ) -\E (X\cdot \nabla U ) , \\
 & = & \E (\mathscr{D}_{\mu} X \nabla K ) +i\mu \E ( \Cor (X_t  ) \cdot \nabla K ) -\E (\alpha U ),\\
 & = & \E (2K -\gamma U +i\mu (\Cor (X_t ) \cdot \nabla K) ) .
\end{array}
\right .
\end{equation}
This concludes the proof.
$\square$\vskip 1mm

Contrary to the strong stochastic virial Theorem, it is not easy to interpret the condition that 
$\di\frac{d}{dt} 
\left ( 
\E (X\cdot \nabla K ) 
\right ) 
=0$ and to ensure that this condition must be fulfilled at equilibrium. 

\begin{remark}
In (\cite{chavanis},V), P.H. Chavanis studies the behavior of $I= m\mathbb{E} (X_t^2 )$ corresponding to the moment of inertia for which he derived a virial type result.
\end{remark}

\section{Some properties of the differential stochastic Newton equation}
\label{properties}

We prove that the stochastic processes which are solutions of the stochastic Newton equation are always gradient diffusions. The complex nature of the stochastic derivative as well as the fact that the stochastic force has real value is fundamental. This result justify an assumption made by E. Nelson in \cite{nelson}. A different argument is given by L. Nottale in (\cite{nottale1},p.149). When the underlying force is dissipative, this argument is no longer true and the solutions of the stochastic Newton equations do not belong to gradient diffusions. This result has strong consequences, as it cancels the possibility to define a natural action functional for the stochastic system which is fundamental in Section \ref{sectionaction}.

\subsection{Conservative forces - Gradient diffusions and the reality condition}

A main feature of the $\alpha$-stochastic Newton equation in the conservative case is that, by construction the potential $U (X_t )$ entering in equation (\ref{dynsto}) is real. This condition gives strong constraints on the stochastic processes solutions of equation (\ref{dynsto}). Indeed, using the composition Lemma we obtain the following constraint on solution of the $\alpha$-stochastic Newton equation:

\begin{lemma}[Reality condition]
\label{realcondition}
A stochastic process $X_t$ is a solution of the $\alpha$-stochastic Newton equation if and only if 
\begin{equation}
(D_+^2 -D_-^2 )[X_t] (1+\alpha ) +(D_- D_+ -D_+ D_- )[X_t] (1-\alpha ) =0 .
\end{equation}
\end{lemma}

\vskip 2mm {\bf Proof}. 
As the potential $U$ is assumed to be real, this implies that the imaginary part of $\mathscr{D}_{\alpha \mu} \circ \mathscr{D}_{\mu} [X_t]$ is real. The result then follows from the composition Lemma.
$\square$\vskip 1mm

This condition reduces to $(D_+^2 -D_-^2 )[X_t] =0$ when $\alpha=1$ corresponding to the differential stochastic Newton equation.\\

When $\alpha=-1$ corresponding to the variational stochastic Newton equation, the reality condition reduces to $(D_- D_+ -D_+ D_- )[X_t] =0$.  

\begin{lemma}
\label{reality}
Let $X_t$ be a solution of the differential stochastic Newton equation. Then, we have $D_+^2 X =D_-^2 X$.
\end{lemma}

\vskip 2mm {\bf Proof}. 
This follows from the fact that in the differential case, we have $\alpha=1$. Lemma \ref{realcondition} gives the result.
$\square$\vskip 1mm

Although this condition does not seem to be fundamental this is precisely the key point as it implies that the drift must be a gradient, a condition that is most of the time assumed in the literature. Here, this follows directly from the structure of the equation. Precisely, we use the main result of \cite{dn1}:

\begin{theorem}[\cite{dn1},Theorem 5] 
\label{darses}
Let $X$ of the form (\ref{proc}), verifying assumption (H), such that $b\in C^2 (\R^d )$ with bounded derivatives and such that for all $t\in (0,T)$ the second order derivatives of $\nabla \ln p_t$ are bounded. Then we have the following equivalence : 
\begin{equation}
D_+^2 X =D_-^2 X \ \mbox{\rm for almost all t}\in (0,T)\ \iff\ b\ \mbox{\rm is a gradient.}
\end{equation}
\end{theorem}

We then obtain the main result of this Section:

\begin{theorem}
\label{gradient}
Any solution of the differential stochastic Newton equation is a gradient diffusion meaning that there exists a potential $W:\R^d \rightarrow \R$ so that the stochastic process is of the form
\begin{equation}
dX_t =\nabla W +\sigma dW_t .
\end{equation}
\end{theorem}

\vskip 2mm {\bf Proof}. 
Let $X_t$ be a solution of the differential stochastic Newton equation, then by Lemma \ref{reality}, we have $D_+^2 X_t = D_-^2 X_t$ so that using Theorem \ref{darses} there exists a certain potential $W:\R^d \rightarrow \R$ so that the stochastic process is of the form $dX_t =\nabla W +\sigma dW_t$.
$\square$\vskip 1mm 

This Lemma allows us to introduce a special function called the {\bf action functional} in Section \ref{sectionaction}.\\

It also proves that the assumption made by E. Nelson in (\cite{nelson2},p.1082, equation (44)) that the real part of $\mathscr{D}_{\mu} X_t$ is a gradient is direct consequence of the structure of the stochastic Newton equation and is then unnecessary. This is used in Section \ref{sectionschro} in order to generalize our previous result \cite{cd} concerning the connection between solutions of the stochastic Newton equation and the Schr\"odinger equation. 

\subsection{Dissipative forces - non existence of a gradient diffusion} 

As already discussed in Section \ref{newtondiff}, considering dissipative forces as a friction or damping term which depends linearly on the speed like $F(x)=-\gamma \dot{x}$ leads by the stochastic embedding to a complex valued force given by $F(X_t )= -\gamma \mathscr{D}_{\mu} X_t$. \\

A main consequence from Theorem \ref{darses} is the following:

\begin{theorem}
Let $F (x,\dot{x})$ be a force such that $F(X ,\mathscr{D} X )$ is not real, then the stochastic processes solution of the stochastic Newton equation are not gradient diffusions. In particular, one can not find a real function $S$ such that the real part of $\mathscr{D} X =\nabla S$. 
\end{theorem}

This Theorem implies directly that for dissipative forces in the stochastic case, the complex speed is never a gradient. \\

Of course, one can artificially recover a gradient by breaking the algebraic structure of the stochastic embedding. Indeed, considering as a generalization of the dissipative force a quantity like 
\begin{equation}
F(X_t , \mathbf{v} ) ,
\end{equation}
where $\mathbf{v} = \di\frac{D_+ X_t + D_- X_t }{2}$, Theorem \ref{darses} applies and we obtain a gradient diffusion. This is for example the strategy followed by P.H. Chavanis in (\cite{chavanis},III.A). Other possibilities can of course be studied. However, there are no reasons a priori to change the way a given equation if embedded over the stochastic processes as long as one is considering dissipative forces.\\

 The fact that $X_t$ is not a gradient in this case is in accordance with the fact that by the Helmholtz Theorem \cite{olver} already in the classical case, no Hamiltonian formulation of the system can be founded \cite{santilli}. In the classical case, one can recover a Hamiltonian formulation by using the fractional calculus \cite{inizan,sina}. One can think to extend such a "fractional" point of view by looking for fractional diffusion processes instead of classical diffusion processes. 

\section{Action function and a stochastic Hamilton-Jacobi equation}
\label{sectionaction}

A useful consequence of Theorem \ref{gradient} is that the speed $V=\mathscr{D}X$ can be written as a gradient of a complex function denoted by $\mathscr{A} (t,X)$, which corresponds in the classic case to the {\bf action functional} \cite{arn}. We prove that the complex action functional satisfies a generalization of the classical Hamilton-Jacobi equation in the stochastic case. Previous work in this direction has been made by F. Guerra and L.M. Morato in \cite{guerramorato1} where a stochastic Hamilton-Jacobi equation is derived in the context of {\bf conservative diffusion} but dealing with a real action functional. This difference will induce strong differences as the complex character of the functional implies the emergence of terms mixing the real and complex part of the functional in the stochastic Hamilton-Jacobi equation.

\subsection{Action functional}

Action functional plays a fundamental role in mechanics (see \cite{arn}, Chap. 9). Using the fact that $\mathscr{D} X_t$ can be written as the gradient of a function, we introduce a natural analogue of the classical action functional in our stochastic setting. As the stochastic derivative is complex, the action functional takes also values in $\C$. This will have a strong impact on the corresponding stochastic Hamilton-Jacobi equation that we derive in the next Section.

\begin{lemma}[Action function]
\label{action}
Let $X$ be a stochastic process solution of the differential stochastic Newton equation (\ref{dynsto}). Then, denoting $\mathscr{V} = \mathscr{D} X$ the complex speed of the process, there exists a function $\mathscr{A} (t,X)$ called the action function such that  
\begin{equation}
\mathscr{V} = \di\frac{\nabla \mathscr{A} (t,X_t )}{m} ,
\end{equation}
with 
\begin{equation}
\mathscr{A} (t,X_t ) = S (t,X_t ) +i\mu R(t,X_t) ,
\end{equation}
where
\begin{equation}
S (t,X)=m W(t,X) -\di\frac{1}{2} m \sigma^2 \ln (p_t ) ,\ \ R(t,X)= \di\frac{1}{2} m \sigma^2 \ln (p_t ) ,
\end{equation}
and the function $W$ is given by Theorem \ref{darses}.
\end{lemma}

\vskip 2mm {\bf Proof}. 
This follows from the properties of the Nelson derivatives. Indeed, we have $D_+ X = V(X_t )$ and 
$D_- X = V(X_t ) -\sigma^2 \nabla \ln (p_t )$. As $X$ is a solution of the differential stochastic Newton equation, we have $V(X_t ) = \nabla W (X_t )$ so that $D_+X$ and $D_-X$ are then gradients. By definition of $\mathscr{D}X$, using the fact that $(D_+ X +D_- X)= 2 V(X_t ) -\sigma^2 \nabla ln (p_t)$ and $(D_+ X -D_- X)= \sigma^2 \nabla \ln (p_t )$ and replacing $V$ by its expression, we obtain the result. 
$\square$

\begin{remark}
In \cite{guerramorato1}, F. Guerra and L.M. Morato do not introduce the previous action functional as they always consider real quantities. However, if a stochastic process $X_t$ is a {\bf critical} diffusion in the sense of Lafferty (see \cite{nelson3},p.439) for a certain action functional then one has that real part of $\mathscr{D} X_t$ is the gradient of a function $S$. We return further on this topic after as they derive a stochastic Hamilton-Jacobi equation for $S$.
\end{remark} 

\subsection{A stochastic Hamilton-Jacobi equation}

The action function $\mathscr{A}$ satisfies in the classical case (i.e. corresponding to $\sigma =0$) a nonlinear first-order partial differential equation called the Hamilton-Jacobi equation (see \cite{arn},Chap. 9, Section 46, p.255). In a similar way, we derive a stochastic version of the Hamilton-Jacobi equation:

\begin{theorem}[A stochastic Hamilton-Jacobi equation]
The action function $\mathscr{A} (t,X)$ satisfies the nonlinear partial differential equation
\begin{equation}
\label{stocjacobi}
\partial_t \mathscr{A} +\di\frac{1}{2m} \left [ \nabla \mathscr{A} \cdot 
\nabla \mathscr{A} \right ]  
+i\mu \di\frac{\sigma^2}{2} \Delta \mathscr{A}=- U .
\end{equation}
\end{theorem}

It must be noted that when the dynamics is not stochastic, meaning that $\sigma=0$, then $\mathscr{A}$ is real and reduces to $S$ as $R=0$ in this case. As a consequence, the stochastic Hamilton-Jacobi equation is equivalent to
\begin{equation}
\label{inter4}
\partial_t S +\di\frac{1}{2m} \left [ \nabla S \cdot 
\nabla S \right ]  
=- U ,
\end{equation}
which corresponds to the classical Hamilton-Jacobi equation introducing the Hamiltonian function 
\begin{equation}
H (p,x,t)=\di\frac{1}{2} p^2 ,
\end{equation}
and rewriting equation (\ref{inter4}) as 
\begin{equation}
\partial_t S +H (\nabla S ,x,t 
)  
=- U .
\end{equation}

\begin{remark}
A stochastic Hamilton-Jacobi equation was derived by F. Guerra and L.M. Morato in the context of stochastic mechanics (see \cite{guerramorato1,guerramorato2}). As their action functional is real, they obtain only a part of our equation. Moreover, their equation does not cover our result due to the fact that the complex nature of $\mathscr{A}$ mixes the real and imaginary part of $\mathscr{A}$ even by restricting our attention to the real part of the stochastic Hamilton-Jacobi equation. We refer to Section \ref{realstocham} for more details.
\end{remark}

\begin{remark}
The stochastic Hamilton-Jacobi equation (\ref{stocjacobi}) has nothing to do with the one introduced in the context of stochastic geometric mechanics in \cite{cami} following the seminal work of J-M. Bismut \cite{bismut}. Indeed, these authors use the Stratonovich stochastic calculus (see \cite{ok}) in order to preserve easily all the geometrical features of the classical Hamiltonian equations. This is not the case here where the It\^o stochastic calculus is used. 
\end{remark}

\vskip 2mm {\bf Proof}. 
Using Lemma \ref{action}, the differential stochastic Newton equation corresponding to $\sigma=1$ can be written as 
\begin{equation}
\mathscr{D}_{\mu} \left [  \nabla \mathscr{A} (t,X_t) \right ] = -\nabla U (X_t ) .
\end{equation}
Using the chain rule formula (\ref{chaine}), we then obtain
\begin{equation}
\label{inter1}
\partial_t \nabla \mathscr{A} +\nabla \left [ \nabla \mathscr{A} \right ] \mathscr{D} X +i\mu \di\frac{\sigma^2}{2} \Delta \nabla \mathscr{A}   =- \nabla U .
\end{equation}
The regularity of $\mathscr{A}$ implies that we have 
\begin{equation} 
\partial_t \nabla \mathscr{A}  =\nabla \partial_t \mathscr{A} .
\end{equation}
Moreover, we have  
\begin{equation}
\left .
\begin{array}{lll}
\nabla \left [ \nabla \mathscr{A} \right ] \mathscr{D} X  & = &  
\nabla \left [ \nabla \mathscr{A} \right ] \di\frac{\nabla \mathscr{A}}{m} , \\
 & = & \di\frac{1}{2m} \nabla \left [ \nabla \mathscr{A} \cdot \nabla \mathscr{A} \right ] ,
\end{array}
\right .
\end{equation}
using the fact that the operator $\nabla$ satisfies the Leibniz relation.\\

We also have 
\begin{equation}
\Delta \nabla \mathscr{A} =\nabla \left [ \Delta \mathscr{A} \right ] .
\end{equation}

As a consequence, equation (\ref{inter1}) can be rewritten as 
\begin{equation}
\label{inter2}
\nabla 
\left [ 
\partial_t \mathscr{A} +\di\frac{1}{2m} \left [ \nabla \mathscr{A} \cdot 
\nabla \mathscr{A} \right ]  
+i\mu \di\frac{\sigma^2}{2} \Delta \mathscr{A}  \right ] =- \nabla U .
\end{equation}

$\square$\vskip 1mm

\section{Emergence of a stochastic induced potential}
\label{sectionemergence}

As $\mathscr{A}$ is complex, we can explicit the real and imaginary part of the nonlinear partial differential system (\ref{stocjacobi}) in term of $S$ and $R$. Precisely, we have:

\begin{theorem}
The real functions $S$ and $R$ satisfy the following system of nonlinear partial differential equations: 
\begin{equation}
\label{syst}
\left .
\begin{array}{l}
\di\frac{\partial S}{\partial t} + \di\frac{1}{2m} (\nabla S )^2 -\mu^2 \left ( \di\frac{1}{2m} (\nabla R )^2 +\di\frac{\sigma^2}{2} \Delta R \right ) + U =0 ,\\
\di\frac{\partial R}{\partial t} + \di\frac{1}{m}(\nabla S ) \cdot (\nabla R ) +\di\frac{\sigma^2}{2} \Delta S =0 .
\end{array}
\right .
\end{equation}
\end{theorem}

It must be noted that the term 
\begin{equation}
-\mu^2 \left ( \di\frac{1}{2m} (\nabla R )^2 +\di\frac{\sigma^2}{2} \Delta R \right ) ,
\end{equation}
which can be interpreted as a modification of the potential $U$ induced by the stochastic character of the motion is intimately related to the complex character of $\mathscr{A}$ as these terms do not exist when $\mathscr{A}$ is real. 

\begin{remark}
In \cite{guerramorato1}, F. Guerra and L.M. Morato obtained the real part of the stochastic Hamiilton-Jacobi equation using a particular notion of critical diffusion called Lafferty critical diffusion by E. Nelson in (\cite{nelson3},p.439) and considering a real action functional mixing in a particular way $D_+ X_t$ and $D_- X_t$. We refer to \cite{guerramorato1} for more details and \cite{nelson3} for a discussion of this work. 
\end{remark}

We can go further by expressing $R$ as a function of $p$.

\subsection{Real part of the stochastic Hamilton-Jacobi equation - emergence of a stochastic induced potential}
\label{realstocham}

A simple calculation leads to the following modified Hamilton-Jacobi equation:

\begin{lemma}[Modified Hamilton-Jacobi equation] 
The first equation of system (\ref{syst}) can be rewritten as 
\begin{equation}
\di\frac{\partial S}{\partial t} + \di\frac{1}{2m} (\nabla S )^2 -m\di\frac{\sigma^4}{2} \di\frac{\Delta (\sqrt{p} )}{\sqrt{p}} +  U =0
\end{equation}
called the modified Hamilton-Jacobi equation.
\end{lemma}

\vskip 2mm {\bf Proof}. 
We denote by $\alpha= m\di\frac{\sigma^2}{2}$. We have $\nabla R =\alpha \di\frac{\nabla p}{p}$ and $\Delta R = -\alpha \left ( \di\frac{\nabla p \cdot \nabla p}{p^2} - \di\frac{\Delta p}{p} \right )$. As a consequence, we obtain 
\begin{equation}
\left .
\begin{array}{lll}
-\di\frac{1}{2} (\nabla R )^2 -m\di\frac{\sigma^2}{2} \Delta R & = &  -\di\frac{\alpha^2}{2} \di\frac{\nabla p \cdot \nabla p}{p^2} +\alpha^2 \left ( \di\frac{\nabla p \cdot \nabla p}{p^2} - \di\frac{\Delta p}{p} \right ) ,\\
  & = & - \alpha^2 \frac{\Delta p}{p} +\di\frac{\alpha^2}{2} \di\frac{\nabla p \cdot \nabla p}{p^2} .
\end{array}
\right .
\end{equation} 
Writing $p$ as $p=(\sqrt{p})^2$ and using the identity $\Delta f^2 = 2f\Delta f + 2\nabla f \cdot \nabla f$, one obtain the identity (with $f=\sqrt{p}$):
\begin{equation}
2 \di\frac{\Delta (\sqrt{p})}{\sqrt{p}} = \di\frac{\Delta p}{p} -\di\frac{1}{2} \di\frac{\nabla p \cdot \nabla p}{p^2}
\end{equation}
which leads to 
\begin{equation}
-\di\frac{1}{2} (\nabla R )^2 -m\di\frac{\sigma^2}{2} \Delta R = -2 \alpha^2 \di\frac{\Delta (\sqrt{p})}{\sqrt{p}} .
\end{equation}
This concludes the proof.
$\square$\vskip 1mm

The main observation is that the stochastic nature of the dynamical system modifies the classical Hamilton-Jacobi equation. The new term can be interpreted as the appearance of a new potential, of pure dynamical origin:

\begin{definition}[Induced stochastic potential] Let $\sigma >0$. We call induced stochastic potential and we denote by $U_{\sigma ,induced}$ the potential defined by  
\begin{equation}
\label{induced}
U_{\sigma, induced} =-m\di\frac{\sigma^4}{2} \di\frac{\Delta (\sqrt{p} )}{\sqrt{p}} .
\end{equation}
\end{definition}

The form of this potential is well known in quantum mechanics (see for example \cite{wheeler}) and is called {\bf Bohm's potential} or {\bf quantum potential} in the literature (see for example the book of J-C. Zambrini (\cite{zambrini2},p.168). It was introduced by D. Bohm in his non-local hidden variable theory for quantum mechanics.

\begin{remark}
This emergent potential is exactly the one obtained by D. Rocha and L. Nottale in \cite{rocha} (see also \cite{nottale2}, equation (12.15) p.521).
\end{remark}

\begin{remark}
One must be careful with the previous result by discussing the dynamical consequences of such an induced stochastic potential. This result is purely stochastic so that the real part of the stochastic Hamilton-Jacobi must not be interpreted as the classical Hamilton Jacobi equation with a potential given by $U+U_{\sigma ,induced}$. Indeed, the dynamical properties of a classical deterministic system whose dynamics is controlled by $U+U_{\sigma ,induced}$ are different from the dynamical properties of a stochastic systems whose underlying dynamics is controlled by the potential $U$. This will be of importance when interpreting such a potential in our application to the flat rotation curves of spiral galaxies in Section \ref{darkgalaxies}.
\end{remark}

\subsection{Imaginary part of the stochastic Hamilton-Jacobi equation - a density equation versus the Fokker-Planck equations}

The second equation corresponds in fact to a classical continuity equation:

\begin{lemma}[A continuity equation]
\label{continuity}
The density $p_t$ of a stochastic process solution of the Newton stochastic differential equation satisfies the following continuity equation
\begin{equation}
m\di\frac{\partial p}{\partial t} +\div (p\, \nabla S ) =0 .
\end{equation}
\end{lemma}

\vskip 2mm {\bf Proof}. 
As $R=m\di\frac{\sigma^2}{2} \ln (p)$, we obtain 
\begin{equation}
\di\frac{\partial R}{\partial t} + \di\frac{1}{m} (\nabla S ) \cdot (\nabla R ) +\di\frac{\sigma^2}{2} \Delta S =
\di\frac{\sigma^2}{2p} \left (m  \partial_t p +\nabla S \cdot \nabla p +p \Delta S \right ) =0.
\end{equation} 
Using the algebraic identity $\div (p V)=p\div (V) +\nabla p \cdot V$ which allows us to rewrite the term $\nabla S \cdot \nabla p +p \Delta S$ as $\div (p\, \nabla S)$, we obtain the result.
$\square$\vskip 2mm

As already remarked by L. Nottale in (\cite{nottale1},p.146), the continuity equation is nothing else than a rewriting of the {\bf Fokker-Planck equation} (see \cite{kp}, Section 2.4 p.68-69) for a stochastic process solution of the differential stochastic Newton equation. Indeed, using results of E. Nelson (\cite{nelson}, equation (3) and (4) p.105), we have that the density $p$ satisfies the classical Fokker-Planck equation
\begin{equation}
\partial_t p = -\mbox{div} ( p D_+ X ) +\di\frac{\sigma^2}{2} \Delta p ,
\end{equation}
and another one corresponding to the {\bf time reversed process} associated to $X$
\begin{equation}
\partial_t p = -\mbox{div} (p D_- X ) -\di\frac{\sigma^2}{2} \Delta p . 
\end{equation}
As a consequence, combining the two equations, we obtain that $p$ must satisfies 
\begin{equation}
2\partial_t p = -\mbox{div} \left ( p \di (D_+ +D_-) X \right ) ,
\end{equation}
which leads to the continuity equation derived in Lemma \ref{continuity}.

\section{The differential stochastic Newton equation as a Schr\"odinger equation}
\label{sectionschro}

In order to apply the previous result in concrete situations and to identify the induced potential $U_{\sigma, induced}$, one needs to have access to the density $p$ or equivalently to the imaginary part of $\mathscr{D}_{\mu} X_t$ of the stochastic process $X_t$ solution of the differential stochastic Newton equation. The main point, already observed by L. Nottale \cite{nottale1}, is that the solutions of the differential stochastic Newton equation are in correspondence, via a simple change of variables using the action functional, with a particular nonlinear partial differential equation which reduces to the linear Schr\"odinger equation in some cases. A complete proof of this correspondence was given in \cite{cd} under the assumption that the real part of $\mathscr{D} X_t$ where $X_t$ is a solution of the stochastic Newton equation, is a gradient. This assumption was also postulated by E. Nelson \cite{nelson} in his derivation of the Schr\"odinger equation in stochastic mechanics. Using the fact that solutions of the stochastic Newton equation are gradient diffusions proved in Section \ref{properties}, we are able to cancel the Nelson assumption and to generalize our previous result \cite{cd}. We also discuss the relation between the usual Madelung transform and the stochastic Hamilton-Jacobi equation. 

\subsection{The stochastic Newton equation as a Schr\"odinger equation}

In the following, we indicate explicitly the dependence of the action $\mathscr{A}(t,x)$ with respect to $\mu$. Following E. Nelson in \cite{nelson}, we introduce the following change of variables: 

\begin{definition}[Wave function] Let $\sigma >0$ and $C$ be a non zero real constant. We call wave function and we denote by $\psi_{\sigma ,\mu}$ the function defined by 
\begin{equation}
\label{wavefunction}
\psi_{\mu} (t,x)=\di e^{i \di\frac{\mathscr{A}_{\mu} (t,x)}{C}} .
\end{equation}
\end{definition}

The constant $C$ plays the role of a normalization constant. We will write $\psi$ instead of $\psi_{\sigma ,\mu}$ in the following. The main point is that this function has only a meaning as long as one considers stochastic processes as in this case $\mathcal{A} (t,x)$ is complex by nature.\\

By definition of $\psi$, we have 
\begin{equation}
\nabla \mathscr{A}_{\mu} = -i C \di\frac{\nabla \psi}{\psi} .
\end{equation}
Since by definition of $\mathscr{A}_{\mu}$, $\mathscr{D}_{\mu} X =\di\frac{\nabla \mathscr{A}_{\mu}}{m}$, we deduce that the differential stochastic Newton equation can be written as
\begin{equation}
iC \mathscr{D}_{\mu} \left [ \di\frac{\nabla \psi}{\psi} \right ] =\nabla U .
\end{equation}
A computation then leads to the following Theorem:

\begin{theorem}[Schr\"odinger formulation]
\label{theoschro} 
If $X$ is a solution of the differential stochastic Newton equation then the wave function $\psi$ satisfies the nonlinear partial differential equation 
\begin{equation}
\label{nonlinearschro}
iC \partial_t \psi -\mu \di\frac{\sigma^2 C}{2} \Delta \psi +\di\frac{C (\mu m\sigma^2 +C)}{2m}  \di\frac{\nabla \psi \cdot \nabla \psi}{\psi} = \psi \cdot U . 
\end{equation}
When $\mu=-1$, we obtain a nonlinear Schr\"odinger like equation 
\begin{equation}
iC \partial_t \psi +\di\frac{\sigma^2 C}{2} \Delta \psi +\di\frac{C (-m\sigma^2 +C)}{2m}  \di\frac{\nabla \psi \cdot \nabla \psi}{\psi} = \psi \cdot U . 
\end{equation}
\end{theorem}

The main point is to observe that the {\bf nonlinearity} induced by the stochastization assumption is of a very {\bf particular form}.

\begin{remark}
{\bf Nonlinear wave mechanics} was initiated by L. De Broglie in \cite{broglie} in order to have a better understanding of the relation between wave and particle (see \cite{broglie2},p.227-231). It can be interesting to explore under different geometric assumptions the class of nonlinearity which arise. 
\end{remark}

The proof of Theorem \ref{theoschro} follows the same line as the corresponding result in \cite{cd} under the gradient assumption of E. Nelson. We provide a complete proof for the convenience of the reader. 

\begin{remark}
Equation (\ref{nonlinearschro}) was rederived by P.H. Chavanis in (\cite{chavanis}, equation (B.1)) without mention of our previous work \cite{cd} where this equation is rigorously proved for the first time. It must be noted that this equation can also be derived in the context of non differentiable deterministic trajectories as in \cite{cr2003} and in the context of multiscale functions as in \cite{cp}.
\end{remark}

\vskip 2mm {\bf Proof}. 
Using the stochastic chain rule formula (\ref{chaine}), we obtain 
\begin{equation}
\label{interschro}
iC \left ( 
\partial_t \left ( \di\frac{\nabla \psi}{\psi} \right ) +\mathscr{D}_{\mu} X \cdot \nabla \left ( 
\di\frac{\nabla \psi}{\psi} \right ) +i\mu \di\frac{\sigma^2}{2} \Delta \left ( \di\frac{\nabla \psi}{\psi} \right ) 
\right )
= -\nabla U .
\end{equation}
We have the equality 
\begin{equation}
\left . 
\begin{array}{lll}
\partial_t \left ( \di\frac{\nabla \psi}{\psi} \right )  & = &  \partial_t (\nabla \psi ) \cdot \di\frac{1}{\psi} +\nabla \psi \cdot \partial_t \left ( \di\frac{1}{\psi} \right )  ,\\
 & = & (\nabla \partial_t \psi ) \cdot \di\frac{1}{\psi} -\nabla \psi \cdot \di\frac{\partial_t \psi}{\psi^2}  ,\\
 & = & (\nabla \partial_t \psi ) \cdot \di\frac{1}{\psi} + \nabla \left ( \di\frac{1}{\psi} \right ) \cdot \partial_t \psi , \\
 & = & \nabla \left ( \di\frac{\partial_t \psi}{\psi} \right ) .
\end{array}
\right .
\end{equation}
Moreover, we have
\begin{equation}
\Delta \left ( \di\frac{\nabla \psi}{\psi} \right ) = \di\frac{\Delta \psi}{\psi} -\di\frac{\nabla \psi \cdot \nabla \psi}{\psi^2} ,
\end{equation}
and 
\begin{equation}
\left .
\begin{array}{lll}
\mathscr{D}_{\mu} X \cdot \nabla \left ( 
\di\frac{\nabla \psi}{\psi} \right ) & = & -i\di\frac{C}{m} \di\frac{\nabla \psi}{\psi} \cdot 
\nabla \left ( \di\frac{\nabla \psi}{\psi} \right ) , \\
 & = & -i\di\frac{C}{m} \nabla \left ( 
 \di\frac{1}{2} \di\frac{\nabla \psi \cdot \nabla \psi}{\psi^2} \right ) .
\end{array}
\right .
\end{equation}
Replacing these expressions in equation (\ref{interschro}), we obtain the result.
$\square$\vskip 1mm

It must be noted that the density $p_t$ of the stochastic process $X_t$ solution of the differential stochastic Newton equation is precisely related to the modulus of the wave function $\psi$ as follows: 

\begin{lemma}[Density]
\label{density}
Let $X_t$ be a solution of the differential stochastic differential equation and $\psi$ its associated wave function then 
\begin{equation}
\ln ( \psi \bar{\psi} ) =-\di\mu \frac{m}{C} \sigma^2 \ln (p_t ) .
\end{equation}
\end{lemma}

\vskip 2mm {\bf Proof}. 
Indeed, by definition of $\psi$, we have 
\begin{equation}
\ln ( \psi \bar{\psi} )  = i \di\frac{1}{C} \left ( \mathscr{A}_{\mu} -\overline{\mathscr{A}_{\mu}} \right ) ,
\end{equation}
which gives using the definition of the action functional in Lemma \ref{action} that 
\begin{equation}
\mathscr{A}_{\mu} -\overline{\mathscr{A}_{\mu}} = 2i \mu R = i \mu \sigma^2 \ln (p_t ) , 
\end{equation}
so that 
\begin{equation}
\ln ( \psi \bar{\psi} ) =-\di\mu \frac{m}{C} \sigma^2 \ln (p_t ) .
\end{equation}
This concludes the proof.
$\square$\vskip 1mm

In order to recover a usual relation between the density of the stochastic process and the modulus of the wave function $\psi$, one needs to choose $K$ in such a way that $-\di\mu \frac{m}{C} \sigma^2=1$. This condition can also be seen as a condition canceling the nonlinearity. Precisely, we have:

\begin{cor}
\label{corschro}
Let $X$ be a solution of the differential stochastic Newton equation, then taking the normalization constant 
\begin{equation}
C= -\mu m\sigma^2 ,
\end{equation}
equation (\ref{nonlinearschro}) reduces to the linear partial differential equation 
\begin{equation}
\label{partialschro}
-i\mu m\sigma^2  \partial_t \psi + m \di\frac{\sigma^4 }{2} \Delta \psi = \psi \cdot U , 
\end{equation}
and 
\begin{equation}
\label{densityprocess}
\left | \psi (t,x) \right |^2 =p_t (x) .
\end{equation}
If moreover, we consider the case $\mu=-1$ then equation (\ref{partialschro}) reduces to the classical linear Schr\"odinger equation
\begin{equation}
\label{schrodinger}
im\sigma^2  \partial_t \psi + m \di\frac{\sigma^4 }{2} \Delta \psi = \psi \cdot U . 
\end{equation}
\end{cor}

The previous result offers the possibility to obtain an explicit expression for the induced potential by solving the linear Schr\"odinger equation for a given potential. A first example in done in the next Section to prove the potentiality of the previous formalism.

\subsection{Stochastic Hamilton-Jacobi equation versus the Madelung transform}

The  {\bf Madelung transform} was introduced in the seminal paper of E. Madelung \cite{made} in order to relate the linear Schr\"odinger equation to a {\bf hydrodynamic type system}. Precisely, let us consider the nonlinear Schr\"odinger equation
\begin{equation}
\label{nls}
i\epsilon \partial_t \psi_{\epsilon} +\di\frac{\epsilon^2}{2} \Delta \psi_{\epsilon} =f(\mid \psi_{\epsilon} \mid ^2 ) \psi_{\epsilon} .
\end{equation}
The Madelung transform is defined as follows:

\begin{definition}
The Madelung transform is the map $\Phi$ which to any pair of functions $\rho :\R^n \times \R \rightarrow \R_{>0}$ and $\theta : \R^2 \times \R \rightarrow \R$ associates a complex-valued function 
\begin{equation}
\Phi :(\rho ,\theta ) \mapsto \psi:= \sqrt{\rho} \di e^{i\theta /\epsilon} .
\end{equation}
\end{definition}

Denoting by $v=\nabla \theta$, we obtain the {\bf hydrodynamic form} of the nonlinear Schr\"odinger equation (\ref{nls}):

\begin{equation}
\left \{ 
\begin{array}{l}
\partial_t v 
v \cdot \nabla v +\nabla f (\rho ) =\di\frac{\epsilon^2}{2} \nabla \left ( \di\frac{\Delta (\sqrt{\rho} )}{\sqrt{\rho}} \right ) , \\
\partial_t \rho +\mbox{\rm div} (\rho v ) =0 .
\end{array}
\right .
\end{equation}

The Madelung transform can be interpreted in our setting as follows:\\

Let $m=1$. The wave function defined in (\ref{wavefunction}) with the constant $C=-\mu \sigma^2$ given in Corollary \ref{corschro} and the condition $\mu=-1$ in order to obtain the linear Schr\"odinger equation, can be written as  
\begin{equation}
\psi = \sqrt{p} \di e^{i S/\sigma^2} .
\end{equation}
As a consequence, we can intrepret the function $\rho$ and $\theta$ of the Madelung transform:\\

\begin{itemize}
\item We have $\rho=p$, $\theta =S$, i.e. the density of the stochastic process $X_t$ solution of the stochastic Newton equation and the real part of the complex speed $\mathscr{D}_{\mu} X_t$ or of the complex action functional $\mathscr{A}_{\mu}$. Taking $\epsilon =\sigma^2$, the hydrodynamic form of the Schr\"odinger equation is nothing else than the stochastic Hamilton-Jacobi equation associated to $\mathscr{A}_{\mu}$.\\

\item The right-hand side of the hydrodynamic form is called the {\bf quantum pressure} and corresponds to the induced stochastic potential $U_{\sigma ,induced}$.
\end{itemize}

A more "physical" discussion of the Madelung transform is given by P.H. Chavanis in (\cite{chavanis},II.C, II.G, II.H and III.B).


\section{Characterization of the set of diffusion processes and the stochastic induced potential}
\label{strategy}

In order to give a constrained theory for the applications of the previous formalism, one needs to identify the diffusion coefficient $\sigma$ from observational data and to give the explicit form of the induced potential. We indicate a strategy, which will be followed in Section \ref{application}, in order to identify $\sigma$ using the previous correspondence the stochastic Newton equation and the Schr\"odinger equation. \\

Let $U$ be a given potential and let $\sigma >0$. The main steps are the following:

\begin{itemize}
\item Write the Schr\"odinger equation (\ref{schrodinger}) and compute the ground state solution $\psi$.

\item Using formula (\ref{densityprocess}), compute the density $p_t (x)$ of the stochastic process $X$.

\item Compute the induced stochastic potential $U_{\sigma ,induced}$ using formula (\ref{induced}).

\item Identify the diffusion coefficient $\sigma$.
\end{itemize}

As one can see, the fact that the Schr\"odinger equation is underlying the dynamics of the stochastic Newton equation allows to bypass the explicit resolution of the stochastic equation by obtaining the density of the stochastic process. 

\section{Application in the Kepler case} 
\label{application}

In this Section, we use the previous formalism to explicit the induced potential when the initial potential is the Kepler potential. The induced potential takes the form of the ad-hoc dark potential used in the literature to explain the flat rotation curves of spiral galaxies. Using the modified Hamilton-Jacobi equation we prove that the real part of the stochastic speed is indeed constant at the equilibrium. We also prove that the expectation of the real part of the stochastic angular momentum is a first integral. This result can be interpreted as the fact that at equilibrium the mean motion takes place in a plane. We then discuss how the previous results are formulated in a polar coordinates systems leading to the fact that if the motion is assumed to be circular then the orthoradial speed is constant.

\subsection{Induced potential in the Keplerian case}

In this Section, we assume that the potential $U$ is given for all $(x,y,z) \in \R^3 \setminus \{ 0\}$ by the Kepler potential 
\begin{equation}
\label{kepler}
U(x,y,z)= \di\frac{-GM m}{r} , 
\end{equation}
where $r=\sqrt{x^2 +y^2 +z^2}$, $M>0$ and $G$ is the universal constant of gravitation. \\

The linear Schr\"odinger equation (\ref{schrodinger}) with a Kepler potential has well known solutions as it corresponds to the Hydrogen atom model. It is well known that the ground state solution of the linear Schr\"odinger equation (\ref{schrodinger}) given in Corollary \ref{corschro} is such that
\begin{equation} 
\left | \psi \right | = \di\frac{C}{m\sigma^4} \di e^{-2r/r_0} ,
\end{equation}
where $C$ is a real constant and $r_0$ is given by 
\begin{equation}
\label{r0}
\di\frac{1}{r_0} = \di\frac{GM}{2\sigma^4} .
\end{equation}
We then have the density $p_t (x)$ of the underlying stochastic process equal to 
\begin{equation}
p_t (x)= \di\frac{C^2}{m^2 \sigma^8} \di e^{-4r/r_0} .
\end{equation}
Using this expression, we have the following Lemma:

\begin{lemma}[Induced potential-Kepler case] Let $\mu=-1$ and $U$ be given by the Kepler potential (\ref{kepler}). Then the induced stochastic potential is given by 
\begin{equation}
U_{\sigma, induced} = -\di\frac{GMm}{r_0} \di\left ( 1-\di\frac{r_0}{r} \right ), 
\end{equation}
with $r_0$ given by (\ref{r0}).
\end{lemma}

\vskip 2mm {\bf Proof}. 
Denoting by $\gamma$ the quantity $\di\frac{C}{m\sigma^4}$, the quantity $\di\sqrt{p_t (x)}$ is given by $\gamma \di e^{-2r/r_0}$. We deduce that 
\begin{equation}
\di\frac{\partial \sqrt{p} }{\partial x} = -2 \di\frac{x}{r_0 r} \sqrt{p} ,
\end{equation}
and 
\begin{equation}
\di\frac{\partial^2 \sqrt{p}}{\partial x^2} = 2\di\frac{1}{r_0 r} \sqrt{p} \di\left ( \di\frac{2x^2}{r_0 r} -1 +\di\frac{x^2}{r^2} \right ) .
\end{equation}
As a consequence, we obtain 
\begin{equation}
\Delta (\sqrt{p} ) =
2\di\frac{1}{r_0 r} \sqrt{p} \di\left ( \di\frac{2r^2}{r_0 r} -3 +\di\frac{r^2}{r^2} \right ) =\di\frac{4}{r_0^2} \sqrt{p} \left ( 1- \di\frac{r_0}{r} \right ) . 
\end{equation}
As a consequence, the induced potential given by $U_{\sigma, induced} =-m\di\frac{\sigma^4}{2} \di\frac{\Delta (\sqrt{p} )}{\sqrt{p}}$ can be explicitly written as 
\begin{equation}
U_{\sigma , induced} =-\di\frac{2m\sigma^4}{r_0^2} \left ( 1- \di\frac{r_0}{r} \right ) . 
\end{equation}
Replacing $\sigma$ by its expression in function of $r_0$ using (\ref{r0}), we obtain the result.
$\square$\vskip 1mm

The form of the induced potential is exactly the dark potential used in order to recover the flat rotation curve of spiral galaxies (see \cite{nottale2} p.652 for a discussion).

\subsection{The flat rotation curves Theorem and the diffusion coefficient}

We are now in position to explore the consequences of the emergence of this extra potential on the dynamics. The main tool is the stochastic Hamilton-Jacobi equation that we have already proved in Section \ref{sectionaction}. 

\begin{theorem}[Flat rotation curves]
\label{flat}
At equilibrium the real part of the speed $\mathscr{D}_{\mu} X$ denoted by $\mathbf{v}$ has a constant norm equal to $v_0$ where $v_0$ is given by 
\begin{equation}
\label{fv0}
v_0^2 =\di\frac{2GM}{r_0} .
\end{equation}
\end{theorem}

\vskip 2mm {\bf Proof}. 
The modified Hamilton-Jacobi equation at equilibrium, i.e. $\partial_t S =0$, gives 
\begin{equation}
\di\frac{1}{2m} (\nabla S )^2  =-U-U_{\sigma, induced} .
\end{equation}
Using the expression of $U_{add}$ just obtained, we have $U+U_{add} = -\di\frac{GMm}{r_0}$. As a consequence, the real part of the speed $\mathscr{D}_{\mu} X$ denoted by $\mathbf{v}$ satisfies 
\begin{equation}
\mathbf{v} \cdot \mathbf{v} =\di\frac{2GM}{r_0} ,
\end{equation}
which is constant.
$\square$\vskip 1mm

A consequence of the previous Theorem is that the diffusion coefficient $\sigma$ can be evaluated as a function of $v_0$ using the definition of $r_0$ given by equation (\ref{r0}):

\begin{lemma}[Diffusion coefficient]
The diffusion coefficient $\sigma$ is given by 
\begin{equation}
\sigma^2 =\di\frac{GM}{v_0} .
\end{equation}
\end{lemma}

The form of $\sigma$ coincides with the one obtained by L. Nottale in (\cite{nottale2}, p.652, (13.153)).\\

The dimension of $\sigma^2$ is $m^2 s^{-1}$ where $m$ stands for meter and $s$ for seconds which is coherent with formula (\ref{r0}) and the fact that $\sigma dW_t$ must be proportional to a distance $m$. Indeed, as $dW_t$ is in $s^{1/2}$ and $\sigma$ in $m s^{-1/2}$ we obtain that $\sigma dW_t$ is in $m$. 

\subsection{A stochastic Noether Theorem and preservation of the stochastic angular momentum}

A useful property of the motion in a central potential is that the motion is restricted to a plane due to the preservation of the angular momentum. It can be interesting to study the preservation of this property under stochastization. Adapting our proof of the stochastic Noether theorem in \cite{cd}, we prove that the real part of the stochastic angular momentum $X\wedge \mathscr{D}_{\mu} X$ is preserved under the motion of the differential stochastic Newton equation.\\

We have the following result:

\begin{theorem}
\label{invariancetheorem}
Let $U$ be a central potential, meaning that $U$ depends only on $\parallel x \parallel$ and $X$ be a solution of the differential stochastic Newton equation, then we have the following identity
\begin{equation}
\label{invariance}
\di\frac{d}{dt} \E \left ( X\wedge \mathscr{D}_{\mu} X \right ) -i\mu \sigma^2  \E 
\left (  \di\frac{\nabla p}{p} \wedge \mathscr{D}_{\mu} X \right ) =0 .
\end{equation}
\end{theorem} 

We postpone the proof at the end of this Section. \\

The quantity $X\wedge \mathscr{D}_{\mu} X$ is the stochastic analogue of the angular momentum $x_t \wedge v_t$. A classical result states that the angular momentum is a constant vector due to the invariance of the Lagrangian (\ref{lagrangefunc}) under the group of rotations. This result is important since it implies that the motion takes place in a plane orthogonal to the angular momentum. This fundamental result extends in the stochastic case for what concerns the real part of the stochastic angular momentum. Precisely, we have:

\begin{lemma}
\label{lemmeangular}
Let $U$ be a central potential and $X$ be a solution of the differential stochastic Newton equation. Let $\mathbf{v}_t$ be the real part of $\mathscr{D}_{\mu} X$, then we have 
\begin{equation}
\di\frac{d}{dt} \E \left  ( X \wedge \mathbf{v} \right ) =0 .
\end{equation}
\end{lemma}

\vskip 2mm {\bf Proof}. 
We have to evaluate the real part of the correcting term $i\mu \sigma^2  \E 
\left (  \di\frac{\nabla p}{p} \wedge \mathscr{D}_{\mu} X \right )$. As $m\mathscr{D}_{\mu} X = m\mathbf{v} +i\mu \nabla R$ by Lemma \ref{action}, the real part is given by 
\begin{equation}
-\sigma^2  \E 
\left (  \di\frac{\nabla p}{p} \wedge \di\frac{\nabla R}{m} \right ) .
\end{equation}
By definition of $\nabla R$ we have $\nabla R =m\di\frac{\sigma^2}{2} \di\frac{\nabla p}{p}$ and the previous quantity reduces to $-2 \E \left ( \nabla R \wedge \nabla R \right ) =0$. As a consequence, taking the real part of equation (\ref{invariance}), we obtain $\di\frac{d}{dt} \E \left  ( X \wedge \mathbf{v} \right ) =0$.
$\square$\vskip 1mm

As a consequence, the real part of the expectation of the stochastic angular momentum 
\begin{equation}
\label{angularstoc}
L_t = \E (X_t \wedge \mathbf{v}_t ) ,
\end{equation}
is a constant vector.\\

The consequence of this result on the motion $X_t$ of the differential stochastic Newton equation are not easy to deduce. However, the situation is very simple at equilibrium:

\begin{theorem}
At equilibrium, the mean motion $E(X_t)$ takes place in a plane orthogonal to the constant vector $L_0$ defined by (\ref{angularstoc}).
\end{theorem} 

\vskip 2mm {\bf Proof}. 
At equilibrium the real part of the speed $\mathscr{D}_{\mu} X$ denoted by $\mathbf{v}_t$ is not random as $\mathbf{v}_t \cdot \mathbf{v}_t$ is constant. As a consequence, we have $\E (X\wedge \mathbf{v} ) =E(X) \wedge \mathbf{v} =L_0$ and the mean motion $E(X_t)$ takes place in a plane orthogonal to $L_0$. 
$\square$\vskip 1mm

\subsubsection{Proof of Theorem \ref{invariancetheorem}}

We consider the classical Lagrangian 
\begin{equation}
\label{newtonlagrangian}
L(x,v)=\di\frac{1}{2} m \parallel v\parallel^2 -\di\nabla U (x) ,
\end{equation} 
where $U$ is assume to be central, i.e. that $U$ depends only on $\parallel x\parallel$.\\ 

We denote by $\{ e_1 ,e_2 ,e_3 \}$ the canonical basis of $\R^3$. Let $\phi_{s,k} : \R^3 \longmapsto \R^3$ be the one parameter family of rotations around the axis $e_k$ for $k=1,2,3$. As an example, for all $X\in \R^3$, we have 
\begin{equation}
\phi_{s,1} (X)=
\left ( 
\begin{array}{ccc}
1 & 0 & 0\\
0 & \cos s & -\sin s\\
0 & \sin s & \cos s 
\end{array}
\right ) 
\cdot X 
.
\end{equation}

We have the following result:

\begin{lemma}
The Lagrangian function $L(X,Z)$, $X\in \R^3$, $Z\in \C^3$ defined by (\ref{newtonlagrangian}) with $U$ a central potential is invariant under the one parameter group of rotations $\phi_{s,k}$ for $k=1,2,3$, i.e. that 
\begin{equation}
L(\phi_{s,k} (X) ,\phi_{s,k} (Z) ) =L(X,Z ) .
\end{equation}
\end{lemma}

\vskip 2mm {\bf Proof}. 
As $\phi_{s,k}$ is an isometry of $\R^3$, we have $\parallel \phi_{s,k} (X)\parallel = \parallel X \parallel$. As $U(x)$ only depends on $\parallel x\parallel$, we then obtain
\begin{equation}
U (\phi_{s,k} (X))=U(X) .
\end{equation}
Extending $\phi_{s,k}$ by linearity to $\C^3$, i.e. $\phi_{s,k} (a+ib )=\phi_{s,k} (a)+i \phi_{s,k} (b)$, we also have $\parallel \phi_{s,k} (Z) \parallel =\parallel Z\parallel$ for $Z\in \C^3$. Indeed, let $Z=a+ib$, $(a,b)\in \R^3 \times \R^3$, then 
\begin{equation}
\parallel \phi_{s,k} (Z) \parallel^2 = \parallel \phi_{s,k} (a) \parallel^2 -\parallel \phi_{s,k} (b)\parallel^2 +2i \phi_{s,k} (a) \cdot \phi_{s,k} (b) .
\end{equation}
As $\phi_{s,k}$ is an isometry of $\R^3$, it preserves the norm of $a$ and $b$ and the scalar product $a\cdot b$. As a consequence, we obtain 
\begin{equation}
\parallel \phi_{s,k} (Z) \parallel^2 = \parallel a \parallel^2 -\parallel b\parallel^2 +2i a \cdot b = 
\parallel Z \parallel^2 .
\end{equation}
As a consequence, we obtain 
\begin{equation}
\di\frac{1}{2} m \parallel \phi_{s,k} (Z) \parallel =\di\frac{1}{2} m \parallel Z \parallel .
\end{equation}
This concludes the proof.
$\square$\vskip 1mm

We then deduce that for all $s\in \R$, the following equality is satisfied
\begin{equation}
L (\phi_{s,k} (X), \phi_{s,k} (\mathscr{D}_{\mu} X ))=L(X,\mathscr{D}_{\mu} X ) .
\end{equation}
By deriving with respect to $s$, we obtain 
\begin{equation}
\di\frac{d}{ds} \left [ 
L (\phi_{s,k} (X), \phi_{s,k} (\mathscr{D}_{\mu} X ))
\right ]
=0 ,
\end{equation}
which can be rewritten as 
\begin{equation}
\label{equa1}
\di\frac{\partial L}{\partial x} ( \phi_{s,k} (X) ,\phi_{s,k} (\mathscr{D}_{\mu} X )) \cdot \di\frac{\partial}{\partial s} ( \phi_{s,k} (X)) +
\di\frac{\partial L}{\partial v} (\phi_{s,k} (X) ,\phi_{s,k} (\mathscr{D}_{\mu} X )) \cdot \di\frac{\partial}{\partial s} (\phi_{s,k} (\mathscr{D}_{\mu} X )) =0 .
\end{equation}
As $\phi_{s,k}$ is linear whose matrix coefficients do not depend on $t$, we have 
\begin{equation}
\phi_{s,k} (\mathscr{D}_{\mu} (X)) = \mathscr{D}_{\mu} (\phi_{s,k} (X)) ,
\end{equation}
and 
\begin{equation}
\di\frac{\partial}{\partial s} \phi_{s,k} (\mathscr{D}_{\mu} (X)) =
\mathscr{D}_{\mu} \left ( 
\di\frac{\partial}{\partial s} \phi_{s,k} (X) \right ) . 
\end{equation}
Moreover for $s=0$, we have $\phi_{0,k} (X)=X$ and equality (\ref{equa1}) in $s=0$ reduces to 
\begin{equation}
\label{equa2}
\di\frac{\partial L}{\partial x} ( X ,\mathscr{D}_{\mu} X ) \cdot \di\frac{\partial}{\partial s} ( \phi_{s,k} (X)) |_{s=0} +
\di\frac{\partial L}{\partial v} (X ,\mathscr{D}_{\mu} X ) \cdot  
\mathscr{D}_{\mu} \left ( \di\frac{\partial}{\partial s} (\phi_{s,k} (X)) |_{s=0} \right ) =0 .
\end{equation}
A simple computation gives 
\begin{equation}
\di\frac{\partial}{\partial s} (\phi_{s,k} (X)) |_{s=0} =e_k \wedge X .
\end{equation}
The differential Newton equation can be rewritten using $L$ as 
\begin{equation}
\mathscr{D}_{\mu} \left ( \di\frac{\partial L}{\partial v} (X,\mathscr{D}_{\mu} (X) ) \right ) =\di\frac{\partial L}{\partial x} (X,\mathscr{D}_{\mu} X ) . 
\end{equation}
Using this equality, we can replace the term $\di\frac{\partial L}{\partial x} ( X ,\mathscr{D}_{\mu} X )$ in equation (\ref{equa2}) by $\mathscr{D}_{\mu} \left ( \di\frac{\partial L}{\partial v} (X,\mathscr{D}_{\mu} (X) ) \right )$. We then have 
\begin{equation}
\label{equa3}
\mathscr{D}_{\mu} \left ( \di\frac{\partial L}{\partial v} (X,\mathscr{D}_{\mu} (X) ) \right )
\cdot \left ( e_k \wedge X \right ) +
\di\frac{\partial L}{\partial v} (X ,\mathscr{D}_{\mu} X ) \cdot  
\mathscr{D}_{\mu} \left ( e_k \wedge X \right ) =0 .
\end{equation}
In order to use the Leibniz formula (\ref{leibnizcor}), we need to explicit the quantity $\Cor (e_k \wedge X)$. A simple computation gives 
\begin{equation} 
\Cor (e_k \wedge X ) =e_k \wedge \Cor (X) ,
\end{equation}
where $\Cor (X) = \sigma^2 \di\frac{\nabla p}{p}$. As a consequence, taking the expectation of equation (\ref{equa3}) leads to 
\begin{equation}
\label{equa4}
\di\frac{d}{dt} \E \left ( 
\di\frac{\partial L}{\partial v} (X,\mathscr{D}_{\mu} (X) ) 
\cdot \left ( e_k \wedge X \right )  \right ) -i\mu \E \left ( 
\di\frac{\partial L}{\partial v} (X ,\mathscr{D}_{\mu} X ) \cdot  
\left ( e_k \wedge \Cor (X) \right ) \right ) =0 .
\end{equation}
As $\di\frac{\partial L}{\partial v} =mv$, we finally have 
\begin{equation}
\label{equa5}
\di\frac{d}{dt} \E \left ( 
\mathscr{D}_{\mu} (X) \cdot \left ( e_k \wedge X \right )  \right ) -i\mu \E\left ( 
\mathscr{D}_{\mu} X \cdot  
\left ( e_k \wedge \Cor (X) \right ) \right ) =0 .
\end{equation}
As we have $u\cdot (v\wedge w)=(u\wedge v) \cdot w$, equation (\ref{equa5}) can be rewritten as 
\begin{equation}
\di\frac{d}{dt} \E \left ( 
e_k \cdot \left (  X \wedge \mathscr{D}_{\mu} (X) \right )  \right ) 
-i\mu \E \left ( 
e_k \cdot \left ( \Cor (X) \wedge \mathscr{D}_{\mu} X  \right ) \right ) =0 ,
\end{equation}
for $k=1,2,3$ which implies that 
\begin{equation}
\di\frac{d}{dt} \E \left ( 
X \wedge \mathscr{D}_{\mu} (X) \right ) 
-i\mu \E \left ( 
\Cor (X) \wedge \mathscr{D}_{\mu} X  \right ) =0 .
\end{equation}
Replacing $\Cor (X)$ by its expression, we obtain the result.

\subsection{Stochastic motion in a plane: radial and orthoradial stochastic speed}

The previous result suggest to study a simplified situation where the stochastic motion is restricted to a plane. In that case, a classical way to describe the motion is to use polar coordinates $(r,\theta )\in \R^+ \times [0,2\pi[$ and to look for the motion $X_t$ as 
\begin{equation}
X_t =r_t e_{\theta_t} , 
\end{equation}
where $e_{\theta_t} = (\cos (\theta_t ) ,\sin (\theta_t ))$. As $X_t$ is a stochastic process, we assume that $r_t$ and $\theta_t$ are diffusion process of the form 
\begin{equation}
dr_t =a_r dt +\sigma_r dW_t,\ \ d\theta_t =a_{\theta} dt +\sigma_{\theta} dW_t ,
\end{equation}
where $a_r$ and $a_{\theta}$ are two functions a priori depending on $r$ and $\theta$ and $\sigma_r$ and $\sigma_{\theta}$ are assumed to be constant in a first approximation. \\

Under these assumptions, one can compute the form of $\mathbf{v}_t$ corresponding to the real part of $\mathscr{D}_{\mu} (X_t )$. Indeed, denoting by $f$ the function defined by  
\begin{equation}
f(r,\theta )= r e_{\theta} ,
\end{equation}
and using the chain rule formula (\ref{chaine}), we have 
\begin{equation}
\left .
\begin{array}{lll}
\mathscr{D}_{\mu} (X_t) & = & \di\mathscr{D}_{\mu} \left [ f(r_t ,\theta_t ) \right ] ,\\
 & = & \mathscr{D}_{\mu} r_t \, e_{\theta_t} +\mathscr{D}_{\mu} \theta_t \, r_t \di e_{\theta_t}^{\perp}
 -\di i\frac{\mu}{2} \sigma_{\theta}^2 e_{\theta_t} +i\mu \sigma_r \sigma_{\theta} e_{\theta_t}^{\perp} ,
\end{array}
\right .
\end{equation}
where $e_{\theta_t}^{\perp} =(-\sin \theta_t , \cos \theta_t )$. \\

The real part of $\mathscr{D}_{\mu} (X_t)$ denoted by $\mathbf{v}_t$ is then given by 
\begin{equation}
\mathbf{v}_t = \mbox{\rm Re} \left [ \mathscr{D}_{\mu} r_t \right ] \, e_{\theta_t} +\mbox{\rm Re} \left [ \mathscr{D}_{\mu} \theta_t \right ] \, r_t \di e_{\theta_t}^{\perp} ,
\end{equation}
where for $z\in \C$  we denoted by $\mbox{\rm Re} \left [z \right ]$ the real part of $z$. \\

As usual, the speed $\mathbf{v}_t$ can be decomposed in two components denoted by $\mathbf{v}_{r,t}$ and 
$\mathbf{v}_{\theta ,t}$ defined by 
\begin{equation}
\mathbf{v}_{r,t} =\mbox{\rm Re} \left [ \mathscr{D}_{\mu} r_t \right ] \, e_{\theta_t} ,\ \ \mbox{\rm and}\ \ 
\mathbf{v}_{\theta ,t} = \mbox{\rm Re} \left [ \mathscr{D}_{\mu} \theta_t \right ] \, r_t \di e_{\theta_t}^{\perp} , 
\end{equation}
corresponding to the radial and orthoradial components of $\mathbf{v}_t$. \\

Using the previous expression, we obtain 
\begin{equation}
X_t \wedge \mathbf{v}_t = r_t^2 \mbox{\rm Re} \left [ \mathscr{D}_{\mu} \theta_t \right ] \mathbf{k} ,
\end{equation}
where $\mathbf{k}$ is the canonical vector $(0,0,1)$. \\

Lemma \ref{lemmeangular} is then equivalent to 
\begin{equation}
\di\frac{d}{dt} \E \left [ r_t^2 \mbox{\rm Re} \left [ \mathscr{D}_{\mu} \theta_t \right ] \right ] =0 ,
\end{equation}
as $\mathbf{k}$ is constant and non zero. \\

The previous result can be used to translate the flat rotation curves Theorem in a polar coordinate systems assuming that the motion is a circular orbit:

\begin{theorem}[Flat rotation curve - polar coordinates]
\label{flatpolar}
At equilibrium, assuming that the motion is a circular orbit of radius $r$, we have $\mathbf{v}_t =v_{\theta}$ and 
\begin{equation}
\E \left [ v_{\theta} \right ] = v_0^2 =\di\frac{L_0^2}{r^2} ,
\end{equation}
where $v_0$ is given by Theorem \ref{flat} and $L_0$ is the constant $L_0 = r^2 \E \left [ \mbox{\rm Re} \left [ \mathscr{D}_{\mu} \theta_t \right ] \right ] $.
\end{theorem}

\vskip 2mm {\bf Proof}. 
We have $\parallel \mathbf{v}_t \parallel^2 = v_r^2 +v_{\theta}^2$. Taking the expectation, we deduce that 
\begin{equation}
\left .
\begin{array}{lll}
\E \parallel \mathbf{v}_t \parallel^2 & = & \E \left [ v_r^2 \right ] + \E \left [ v_{\theta}^2 \right ] , \\
 & = & \E \left [  \mbox{\rm Re} \left [ \mathscr{D}_{\mu} r_t \right ] ^2 \right ] + 
 +\E \left [ r_t^2 \mbox{\rm Re} \left [ \mathscr{D}_{\mu} \theta_t \right ] ^2 \right ] .
\end{array}
\right .
\end{equation}
As the motion is a circular orbit of radius $r$, we have $r_t =r$ and $\mathscr{D}_{\mu} r_t =0$ so that \begin{equation}
\E \parallel \mathbf{v}_t \parallel^2 =\E \left [ v_{\theta}^2 \right ] .
\end{equation}
Moreover, we have $\E \left [ r_t^2 \mbox{\rm Re} \left [ \mathscr{D}_{\mu} \theta_t \right ] \right ] =r^2  \left [\mbox{\rm Re} \left [ \mathscr{D}_{\mu} \theta_t \right ] \right ]$ which is constant equal to $L_0$ by assumption. Then $\E \parallel \mathbf{v}_{\theta} \parallel^2 = r^2 \E \left [ \mbox{\rm Re} \left [ \mathscr{D}_{\mu} \theta_t \right ] ^2 \right ]=L_0^2 /r^2$. 

Theorem \ref{flat} implies that $\E \parallel \mathbf{v}_t \parallel^2 =v_0^2$. This concludes the proof.
$\square$\vskip 1mm

The previous result indicates that in some particular cases, one can directly connect the flat rotation curve Theorem with the fact that the mean of the real part of the orthoradial speed is constant. \\

An interesting consequence of Theorem \ref{flatpolar} and Theorem \ref{flat} using the explicit expression of $v_0$ depending on $r_0$, is the following formula for $r$:
\begin{equation}
r^2 = \di\frac{L_0^2 r_0}{2GM} . 
\end{equation}

\section{The "dark matter" problem and flat rotation curves of spiral galaxies}
\label{darkgalaxies}

As the induced potential in the Kepler case corresponds to the usual ad-hoc dark potential used in astrophysics, one is leaded to discuss the applications of the previous formalism in the setting of the dynamics of galaxies. We follow here the arguments given by D. Rocha and L. Nottale in \cite{rocha} and (\cite{nottale2},Section 13.8.2 p.652-654).

\subsection{Structural assumptions on the galaxy and the dark matter problem}

We consider an isolated spiral galaxy which is already formed. A spiral galaxy is always decomposed in three components:\\

\begin{itemize}
\item The central bulb which looks approximately like a sphere containing gazes and stars in a homogeneous way. 

\item The disk which is itself decomposed in two components: a fine disk which is dense and a rough one. 

\item The pair bulb-disk is contained in a halo of stars which is more or less a sphere but of low density.
\end{itemize}

\begin{center}
\begin{figure}[h!]
\centerline{\includegraphics[scale=0.5]{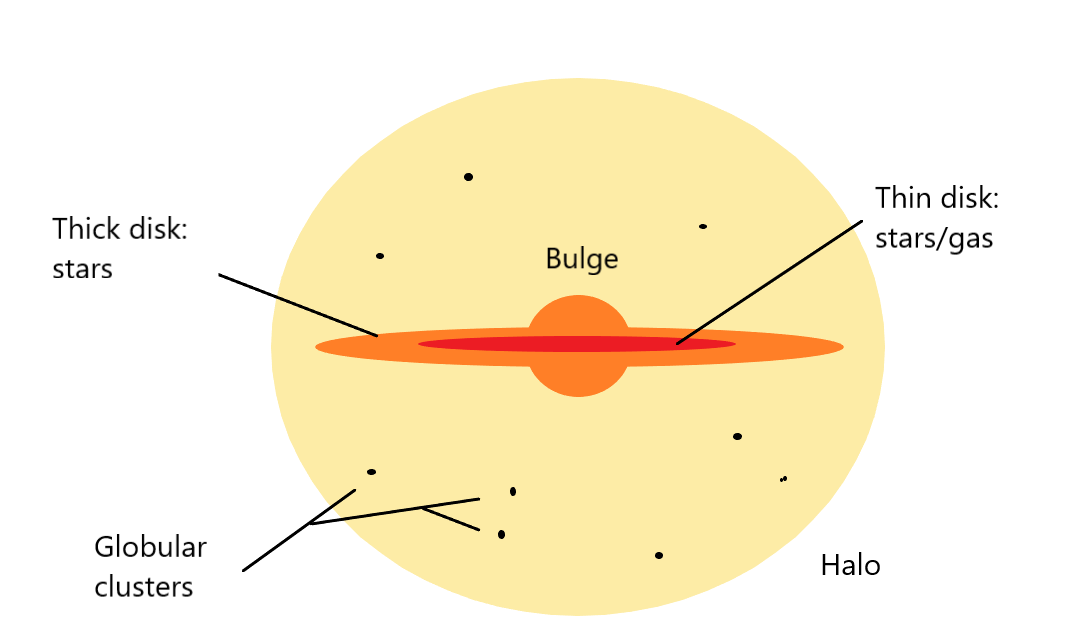}}
\caption{The Milky way.}
\end{figure}
\end{center}

These structures can be observed by all possible methods and give the total mass of the galaxy or in order to be precise the mass which can be detected by our actual observational means. Using these data, classical physics tell us that the stars which are meanly in the disk would have Keplerian orbits around the bulb. In this case, their speed of rotation would decrease as a function of the distance to the center of the galaxy. This is not the case and observational data show that the speed of rotation is more or less constant leading to a paradox. 

\begin{center}
\begin{figure}[h!]
\centerline{\includegraphics[scale=0.8]{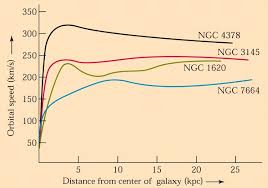}}
\caption{Some examples of rotation curves}
\end{figure}
\end{center}

Two main theories have been designed to solve this problem:

\begin{itemize}
\item The first and most common proposition is the existence of a huge amount of invisible matter with very specific properties called the "dark matter". We refer to \cite{persic} for more details.

\item The second one is a proposal of M. Milgrom \cite{milg} called MOND for MOdified Newton Dynamics consisting in changing the classical Newton law of dynamics asserting a linear relation between the force and acceleration to a non-linear one when the acceleration is weak with respect to a constant fixed by the theory.  
\end{itemize}

A consequence of both approaches is to lead to add a new potential to the classical initial potential. This two point of view have induced a great number of works and we refer to (\cite{blanchet},\cite{beken},\cite{mil2},\cite{sanders}) for more details and recent works. 

\subsection{Brownian diffusion and scale relativity}

In order to apply our results, we need to justify the assumption that we consider the modification of the Newton equation on Brownian diffusion. Such a discussion is provided by S. Albeverio and al. in (\cite{albeverio1},p.366) in the context of the dynamics of a protoplanetary nebulae (see also \cite{cr}, Section 5.1). In our setting, the consideration of the Newton equation in a stochastic setting can be justified in essentially two different points of view:\\

\begin{itemize}
\item In scale-relativity, the fact that space has a fractal structure at small scales implies via the scale-relativity principle that a fractal structure of space also emerges at large scale (see \cite{nottale2},p.559). Assuming that the structure of space outside the bulb becomes fractal with a fractal dimension equal to $2$, one can model the effects of the fractals structure using stochastic processes and in particular Brownian diffusion.  

\item Another point of view, already discuss by L. Nottale in \cite{nottale3} (see also \cite{nottale2}, p.559) is to consider that the underlying dynamics is chaotic, allowing a description of the long-term dynamics using stochastic processes.  
\end{itemize}

In both cases, the use of stochastic processes and in particular Brownian diffusion is a first possible model.\\ 

An essential assumption in order to consider Brownian diffusion with constant coefficient is to assume that the underlying stochasticity is:
\begin{itemize}
\item Isotropic
\item Homogeneous
\end{itemize}
With respect to a general diffusion coefficient $\sigma (t,x)$, the isotropy condition implies that $\sigma$ depends only on $t$ and homogeneity implies that $\sigma (t)$ is constant in time. \\

Of course, the assumption that the diffusion coefficient is constant is only a first approximation and can be weakened in future explorations. In particular, it seems that isotropy is a too strong condition for some applications. 

\subsection{Numerical estimates for a "typical" galaxy}

Assuming that the previous assumptions are valid as a first approximation, one can give some numerical estimates for the diffusion coefficient as long as one can predict the value of the mass $M$, $r_0$ and $v_0$ for a given galaxy.\\

We first introduce some notations and Units:\\ 

\begin{itemize}
\item The mass $M$ of a galaxy is given in Solar mass denoted by $M_{\odot}$ with $M_{\odot} = 1.98 \, 10^{30}$ $kg$. 

\item The size of a galaxy is given in Parsec denoted $Pc$ and equal to $3.086 \, 10^{13}\ km$. 

\item The gravitational constant $G$ is equal to $4.3\, 10^{-6}$ $Kpc\, km^2  s^{-2} M_{\odot}^{-1}$.
\end{itemize}

A typical spiral galaxy has the following properties:\\

\begin{itemize}
\item an extension/size denoted by $l_0$ between $2$ and $100$ $Kpc$. 

\item The mass is of the form $10^p\ M_{\odot}$ where $p$ is between $8$ and $12$. 
\end{itemize}

A typical rotation curve for a spiral galaxy is characterized by two quantities:\\

\begin{itemize}
\item A distance $r_{0,obs}$ at which the rotation curve begins to be flat.

\item A speed $v_{0,obs}$ corresponding to the averaged speed in the flat part of the rotation curve.
\end{itemize}

This two sentences must be taken with some care as the rotation curve is fluctuating around a given value. \\

The typical value for $v_{0,obs}$ is $144$ $km.s^{-1}$ as taken from the Persic-Salucci catalog \cite{cata}. A discussion of this value is given in (\cite{nottale2},p.653-654).

\begin{center}
\begin{figure}[h!]
\centerline{\includegraphics[scale=0.5]{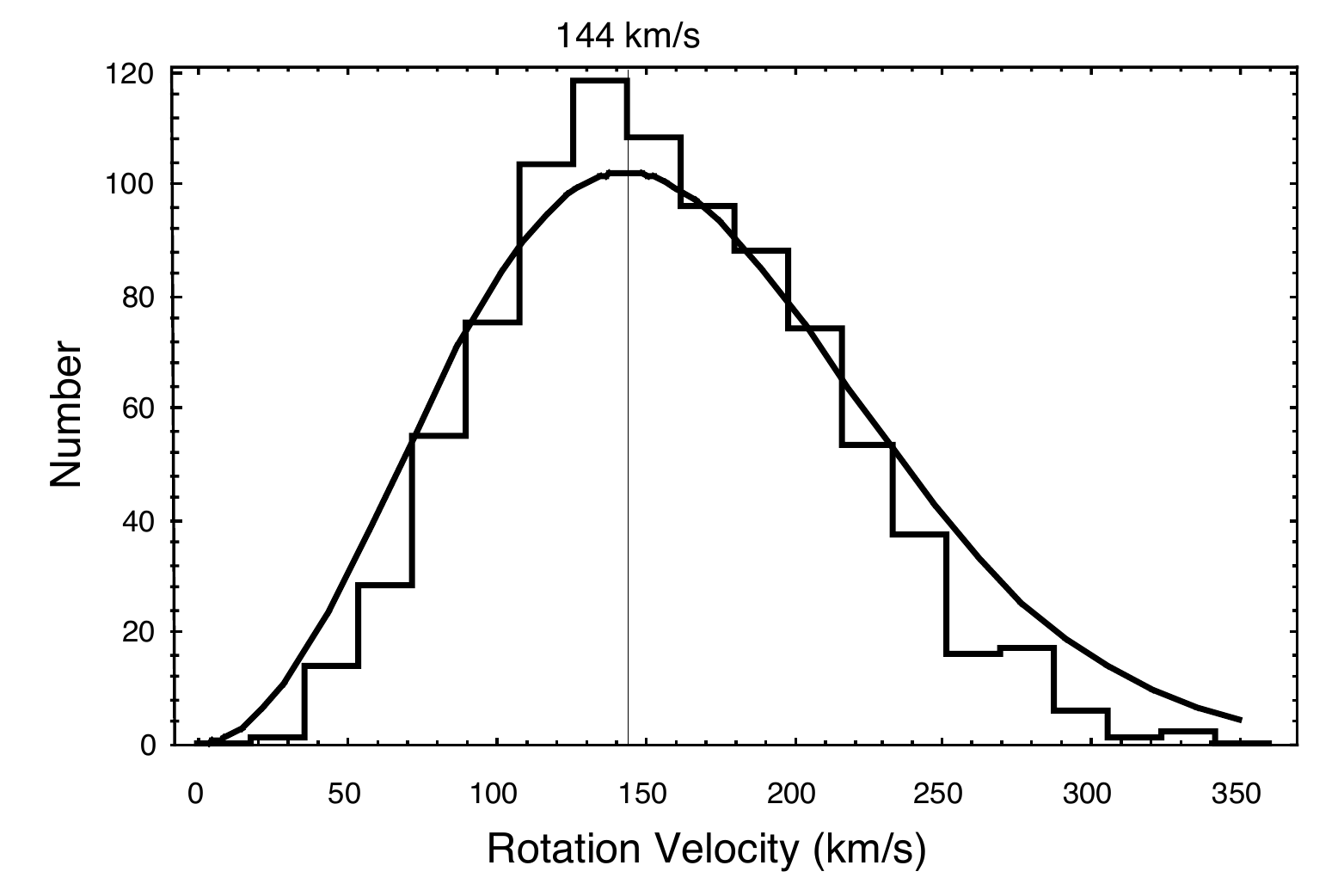}}
\caption{Distribution of rotation velocities of spiral galaxies.}
\end{figure}
\end{center}

If we use our formula relating $r_0$ and $v_0$ taking for $v_0$ the value $v_{0,obs}$, we obtain an estimation denoted $r_{0,estimated}$ of size
\begin{equation}
r_{0,estimated}  = 4.1\, 10^{p-10}\ Kpc .
\end{equation}
As $p$ goes from $8$ to $12$, we then have $r_{0,estimated}$ running from $41\ Pc$ to $410\ Kpc$.\\

For our Galaxy however since the {\bf Visible Mass} is about $8 10^{10} Ms$ and $v_{0,obs}$  is $220$ $km.s^{-1}$ the  $r_{0,estimated}$ is  around $8\ Kpc$ which corresponds approximately to the distance at which the rotation curve begins to be flat as it can be seen with the observations reported in  figure 4.

\begin{center}
\begin{figure}[h!]
\centerline{\includegraphics[scale=0.6]{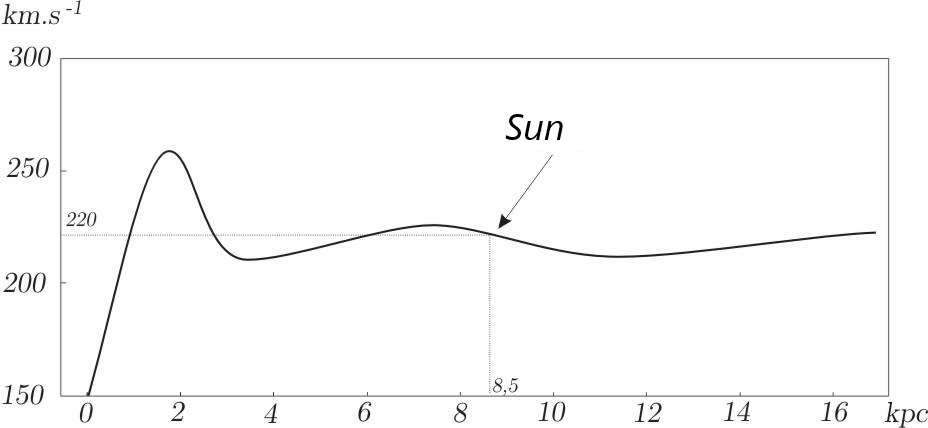}}
\caption{The rotation curve for the Milky Way.}
\end{figure}
\end{center}

The diffusion coefficient $\sigma$ can also be estimated using this value of $v_{0,obs}$ and we obtain
\begin{equation}
\sigma_{estimated}^2 =205 \, 10^{p+3}\ Kpc^2\, s^{-1} ,
\end{equation}
which gives a $\sigma_{estimated}$ of order $14.3\ \di 10^{\frac{p+3}{2}} \ Kpc\, s^{-1/2}$.\\

All the previous computations are based on a very crude approximation on the distribution of matter. However, one can already derive some 
physical conclusion at  this level:\\

The induced potential  applied  to the ordinary visible matter  of a given spiral galaxy (like  ours)  is thus doing the work which is  usually attributed to ad hoc ``dark matter’’. Note that here this potential  is repulsive  for $r< r_{0}$, but it is  attractive  if  not. 

\begin{center}
\begin{figure}[h!]
\centerline{\includegraphics[scale=0.8]{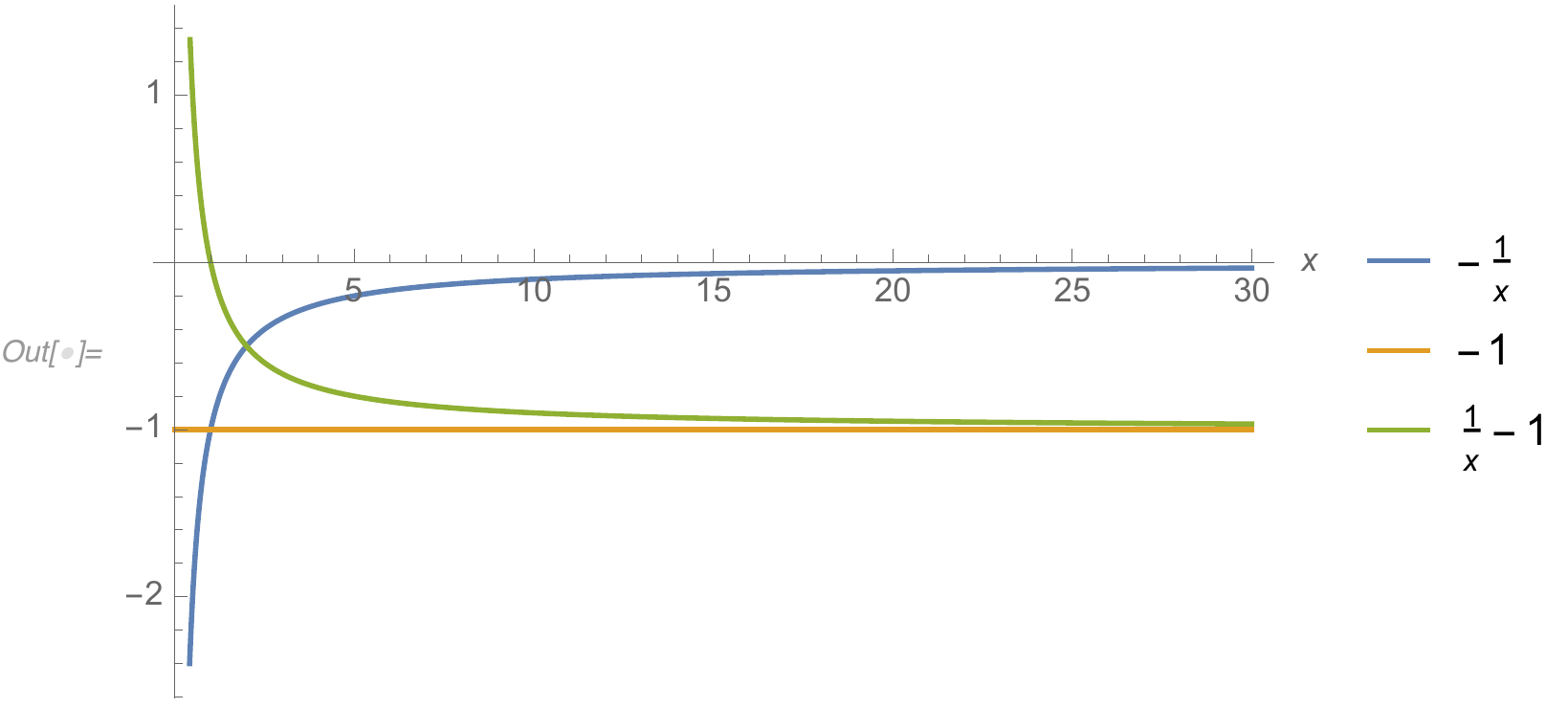}}
\caption{For $r_{0}=1$,  curves for the $3$ potentials  (normalized to $GM/ r_{0}$): $U_{tot} =U_{Kepler} +U_{induced}=-1$ (in orange), $U_{Kepler}$  (in blue) and $U_{induced}$ (in  green).}
\end{figure}
\end{center}

Indeed if one looks at its expression it reaches closely  the required saturation value  as soon as $r$ overtakes a few times the associated $r_{0,estimated}$ value. For example in our Milky Way depending on the chosen model for dark matter and its spatial  extension its mass necessary to explain  the  $v_{0,obs}$  is comprised in the range  $1.5-2.3$ 1$0^{12}$ $M_{\odot}$, obviously not using our relation between  $v_{0}$ and $r_{0}$. Our result confirms indeed the previous work of L. Nottale and D. Rocha in \cite{rocha}. 

\section{Perspectives}
\label{sectionpers}

The previous illustration allows us to put in evidence the main difficulties in order to apply the previous theory:\\

\begin{itemize}
\item Compute the diffusion parameter $\sigma$ and if possible, relates it to observational data.

\item Identify the potential $U$ entering in the Newton equation.

\item Compute an explicit expression for solutions of the corresponding Shcr\"odinger equation.
\end{itemize}

Depending on the physical problem one is studying, this questions can be very complicated.\\

However, the previous results indicates that the strategy proposed in \cite{not} and \cite{rocha} can be rigorously founded in the framework of the stochastic embedding formalism of \cite{cd} based on E. Nelson stochastic derivatives \cite{nelson}.\\

We can generalize the previous work in various directions.\\

A first one is to consider a more general class of Brownian diffusion considering non constant diffusion, i.e. a $\sigma$ depending on $x$. This assumption will lead us to consider a distribution of mass with a density $\rho (x)$ and then a mass $M(x)$ in the computation via the Poisson equation. \\ 

Another possibility is to use the extension of the stochastic derivative to cover fractional diffusion processes as defined in \cite{darses2}.\\

A third one is to explore more general family of potentials and also to take care of dissipative effects in the initial Newton equation. \\

Last, but not least, one can extend the framework of the stochastic embedding to equations on Riemannian manifolds in order to cover a stochastic generalization of relativistic effects. A starting point can be to use the work of T. G. Dankel \cite{dankel}, D. Dohrn and F. Guerra \cite{dohrn} and T. Zastawniak \cite{zas} where the stochastic mechanics of E. Nelson was extended on Riemannian manifolds.

\section*{Acknowledgements}

The authors thanks J-C. Zambrini for references and a guided tour around the history of stochastic mechanics and P.H. Chavanis for comments and references. J. Cresson thanks the GDR "G\'eom\'etrie diff\'erentielle et m\'ecanique" and in particular A. Hamdouni for discussions and the GDR TraG ("Trajectoire Rugueuses") as well as C. Darrigan for support and discussion.

\end{document}